\documentclass[11pt]{article}

\usepackage{booktabs}
\usepackage[margin=15pt,font=small,labelfont={bf,sf},justification=justified]{caption}
\usepackage[superscript,biblabel]{cite}
\usepackage{color}
\usepackage{float}
\usepackage{graphicx}
\usepackage{lscape}
\usepackage[bbgreekl]{mathbbol}
\usepackage{mathtools}
\usepackage{multibib}
\usepackage{multirow}
\usepackage{tabularx}
\usepackage{titleref}
\usepackage{url}
\usepackage{lipsum}
\usepackage{bm}
\usepackage{tikz}
\usepackage{subfig}
\usepackage{caption}

\usepackage{titling}
\usepackage{authblk} 
\pretitle{\begin{center}\bfseries\sffamily\LARGE}
\posttitle{\end{center}}
\usepackage{abstract}

\usepackage{sectsty} 
\allsectionsfont{\sffamily}

\usepackage[left=1in,right=1in,top=1in,bottom=1in,includefoot,heightrounded]{geometry}

\title{An efficient and accurate semi-implicit time integration scheme for dynamics in nearly- and fully-incompressible hyperelastic solids}
\author[1]{Edward M. Terrell}
\author[1,2,3,*]{Boyce E. Griffith}
\affil[1]{Department of Mathematics, University of North Carolina at Chapel Hill, Chapel Hill, NC, USA}
\affil[2]{Department of Biomedical Engineering, University of North Carolina at Chapel Hill, Chapel Hill, NC, USA}
\affil[3]{McAllister Heart Institute, University of North Carolina at Chapel Hill, Chapel Hill, NC, USA}
\affil[*]{Corresponding author, boyceg@email.unc.edu}
\date{}

\providecommand{\keywords}[1]{\textbf{{Keywords---}} #1}

\begin{document}
\maketitle

\begin{abstract}
\noindent The choice of numerical integrator in approximating solutions to dynamic partial differential equations depends on the smallest time-scale of the problem at hand. Large-scale deformations in elastic solids contain both shear waves and bulk waves, the latter of which can travel infinitely fast in incompressible materials. Explicit schemes, which are favored for their efficiency in resolving low-speed dynamics, are bound by time step size restrictions that inversely scale with the fastest wave speed. Implicit schemes can enable larger time step sizes regardless of the wave speeds present, though they are much more computationally expensive. Semi-implicit methods, which are more stable than explicit methods and more efficient than implicit methods, are emerging in the literature, though their applicability to nonlinear elasticity is not extensively studied. In this research, we develop and investigate the functionality of two time integration schemes for the resolution of large-scale dynamics in nearly- and fully-incompressible hyperelastic solids: a Modified Semi-implicit Backward Differentiation Formula integrator (MSBDF2) and a forward Euler / Semi-implicit Backward Differentiation Formula Runge-Kutta integrator (FEBDF2). We prove and empirically verify second order accuracy for both schemes. The stability properties of both methods are derived and numerically verified. We find FEBDF2 has a maximum time step size that inversely scales with the shear wave speed and is unaffected by the bulk wave speed---the desired stability property of a semi-implicit scheme. Finally, we empirically determine that semi-implicit schemes struggle to preserve volume globally when using nonlinear incompressibility conditions, even under temporal and spatial refinement. 
\end{abstract}

\keywords{nonlinear elasticity, incompressible materials, finite elements, semi-implicit }

\newenvironment{itquote}
  {\begin{quote}\itshape}
  {\end{quote}\ignorespacesafterend}

\section{Introduction}
\label{sec:sample1}
\subsection{Background}
Efficient and accurate numerical methods for simulating dynamic, large-scale deformations in nearly-incompressible or incompressible solids remain scarce, despite the ubiquity of such materials in various scientific fields. Many biomaterials, such as skin and muscle, are mainly comprised of water, and are often approximated as incompressible \cite{Ita2022, Mathur2001}. Various polymers, most notably rubber, also exhibit behavior that can be generalized as incompressible \cite{Schaefer2010}.  Fields that heavily depend on simulations involving such incompressible materials, like physiological modeling or biomedical engineering, would benefit greatly from having access to efficient numerical methods that can resolve incompressible dynamics in solids. 

 Numerous finite element schemes exist for problems in finite elasticity. Because of their utility in structural engineering, particular focus has been applied to developing methods that solve for the equilibrium states of solids under static stresses \cite{Simo1984,Coombs2010,Brink1996, Hughes2000,Karabelas2020}. A broadly-used subclass of static methods, known as quasi-static methods, was developed \cite{Argyris1981, Teran2005} to model elastic deformations under gradual dynamic loads. While extremely effective for analyzing deformations in regimes where inertial effects are negligible, these models cannot resolve finer time-scale phenomena like elastic wave propagation \cite{Quintal2011,Arriaga2007}. Evaluating static and quasi-static formulations for solid mechanics reduces to solving nonlinear elliptic PDEs.

 Because of the inherent numerical stiffness associated with resolving nonlinear dynamics in nearly- and fully-incompressible solids, explicit time-stepping methods are often a poor choice due to severe restrictions on time step size \cite{Leveque2007}. Implicit methods are capable of capturing the relevant dynamics without crippling time step size restrictions, but they depend on nonlinear solvers \cite{Rossi2016, Kadapa2021}. This highlights the potential utility of an intermediate form of integrator: a semi-implicit method that is more stable than an explicit method and more efficient than an implicit method.

Extant semi-implicit time integration schemes for solving dynamic problems in nonlinear incompressible elasticity come in various forms, though they all involve the introduction of a Lagrange-multiplier pressure field that is used to enforce incompressibility. Lahiri \textit{et al.} formulated a fractional-step scheme reminiscent of a St\o{}rmer-Verlet or leap-frog integrator \cite{Lahiri2005,Verlet1967}. Gil \textit{et al.} introduced a three-field Petrov-Galerkin scheme coupled with a total variation diminishing Runge-Kutta method \cite{Gil2014,Gottliev2003}. This work also utilizes a variational stabilization technique similar to the dynamic variational multiscale stabilization introduced by Scovazzi \textit{et al.} \cite{Scovazzi2016, Zeng2017} that permits the usage of linear finite elements for all relevant finite element spaces. Kadapa introduced a semi-implicit scheme based on an explicit form of the Newmark-$\beta$ method \cite{Kadapa2021,Chung1994}. These integrators can use time step sizes that scale with the shear wave speed of the modeled material, enabling them to step over dynamics associated with volumetric deformations. However, none of these methods have been demonstrated to achieve the degree of temporal accuracy claimed, nor have their abilities to preserve volume been tested for large deformations. 

 This paper introduces two semi-implicit integration schemes based on the second-order semi-implicit backward differentiation formula (SBDF2), which resembles the standard second-order backward differentiation formula \cite{Britz1997} except that it treats non-stiff force terms explicitly \cite{Chow2021}. Second-order Adams-Bashforth extrapolations are used to push these non-stiff terms into the forward time step \cite{Bashforth1883}. While SBDF2 is well-studied\cite{Keita2024,Ascher1995,Varah1980}, its application to hyperbolic problems remains relatively unexplored, and it has not been applied in the context of nonlinear incompressible elastodynamics. The proposed modifications to SBDF2 are unique in that they leverage explicit extrapolations within the volumetric energy to ensure that all relevant forcing terms are calculated at the future time step, maintaining second-order accuracy while minimizing computational cost.

Semi-implicit schemes for nonlinear elastodynamics are inspired by methods commonly used in computational fluid dynamics---the most obvious being projection methods of the form popularized by Chorin \cite{Chorin1968}. Projection methods require precise enforcement of boundary conditions for the pressure field, which is difficult for general problems in nonlinear elasticity. The schemes proposed in this work distance themselves from the class of projection methods by never explicitly enforcing pressure boundary conditions; they are instead enforced simultaneously with the boundary conditions for the velocity.

The key contributions of this paper involve the validation of semi-implicit schemes in their application to problems in nonlinear incompressible elastodynamics. We used both local-truncation error analysis and empirical convergence studies to validate that the proposed schemes are second order accurate for both compressible and incompressible materials. Furthermore, we determined via linear stability analysis that one of the proposed schemes has a time step size that inversely scales with the shear wave speed of the material and is unaffected by the bulk wave speed of the material, which is the standard for semi-implicit schemes. The other proposed scheme is found to be unconditionally unstable despite being derived from a well-established integrator. Lastly, both proposed schemes are found to struggle with preserving volume in incompressible dynamics if a nonlinear formulation for the volumetric energy is used; using a linear incompressiblity condition enables convergence to volume preservation under temporal and spatial refinement.

 \section{Preliminaries}
\subsection{Continuum mechanics}\label{Mechanics}
    The aim of the algorithms presented in this manuscript is to efficiently model the dynamic deformation of a domain or body, $\Omega_t\in\mathbb{R}^{d}\text{ for }d=2$ or $3$  from time $t=0$ to $T$. The state of $\Omega_t$ at these time points is denoted by $\Omega_0$ and $\Omega_T$ respectively. The \textit{reference} coordinates used to describe position in $\Omega_0$ are $\mathbf{X}$, and the corresponding position of the point $\mathbf{X}$ in \textit{deformed} coordinates after a deformation is described by $\mathbf{x} = \bm{\chi}(\mathbf{X},t)$, in which $\bm{\chi}(\textbf{X},t)$ is an invertible mapping such that $\bm{\chi}:\Omega_0\rightarrow\Omega_t$. The displacement vector, $\textbf{U}=\bm{\chi}(\mathbf{X},t)-\textbf{X}$, is the primary variable for many numerical methods in finite elasticity  \cite{Simo1984}, including the method provided in this paper. 

    The behavior of the displacement field, $\textbf{U}$, is described by several key variables. The deformation gradient tensor, $$\mathbb{F}\equiv\frac{\partial \bm{\chi}}{\partial \textbf{X}}=\nabla_0\mathbf{x}=\mathbb{I}+\nabla_0\textbf{U},$$ is a two-point tensor that describes local deformations near a given reference coordinate, $\textbf{X}$. The Jacobian of deformation, $J = \det \mathbb{F},$
    describes the volumetric change of a neighborhood around material coordinate $\textbf{X}$ from time $0$ to time $t$. A representation for the isochoric deformations encoded by $\mathbb{F}$ is described by $\bar{\mathbb{F}}=J^{-1/3}\mathbb{F}$. The cofactor of $\mathbb{F},$ $\mathbb{H}=J\mathbb{F}^{-T},$
    describes the change in the area of a plane local to a point $\mathbf{X}$ from time $0$ to time $t$ via Nanson's formula \cite{Tadmor2012,Capecchi2007}.

Utilizing a particular hyperelastic material model provides an explicit formulation for a strain energy density function \cite{Kim2011,Pence2015}: $$\mathcal{W}(\mathbb{F},J)=\mathcal{W}^{\text{dev}}(\mathbb{F},J;\mu)+\mathcal{W}^{\text{vol}}(J;\kappa)$$ Given a form for $\mathcal{W}$, the first Piola-Kirchhoff (PK1) stress tensor is calculated as
\begin{equation}\label{PK1}
    \mathbb{P}=\frac{\partial \mathcal{W}}{\partial \mathbb{F}}=\frac{\partial \mathcal{W}^{\text{dev}}}{\partial \mathbb{F}}+\frac{\partial \mathcal{W}^{\text{vol}}}{\partial J}\mathbb{H}.
\end{equation}
The PK1 stress tensor is a two-point tensor that relates forces in the deformed configuration to areas in the undeformed configuration \cite{Tadmor2012,Capecchi2007}. This work employs an isotropic modified Neo-Hookean material model, defined as 
\begin{equation}
    \mathcal{W}^{\text{dev}}(\bar{\mathbb{F}}) := \frac{\mu}{2}(\bar{\mathbb{F}}:\bar{\mathbb{F}}-3) =  \frac{\mu}{2}(J^{-2/3}\mathbb{F}:\mathbb{F}-3). \label{Wdev}
\end{equation}
Defining the volumetric energy using $\bar{\mathbb{F}}$ instead of $\mathbb{F}$ guarantees that the deviatoric energy is computed exclusively from isochoric deformations. More complicated phenomenological material models, like the Mooney-Rivlin model \cite{Mooney} or the Ogden model \cite{Ogden1972}, are capable of more accurately modeling large deformations in rubber-like materials, and models like that of Holzapfel, Gasser, and Ogden \cite{Holzapfel2000} are capable of capturing anisotropic dynamics. A Neo-Hookean material model was used in this work purely for its simplicity; the derivation of the upcoming numerical schemes are not dependent on the formulation used for internal energy. The numerical schemes presented in this work can be easily adapted to accommodate more complex constitutive models. 

We consider two forms of volumetric energy: \begin{equation}
   \mathcal{W}^{\text{vol,q}}(J;\kappa):= \frac{\kappa}{2}(J-1)^2 \qquad \mathcal{W}^{\text{vol,L}}(J;\kappa):= \kappa(J\ln J -J + 1)   . \label{Wvol}
\end{equation}
$\mathcal{W}^{\text{vol,q}}$ is the most commonly utilized formulation for the volumetric energy due to its simplicity. $\mathcal{W}^{\text{vol,L}}$ was first proposed by Liu \textit{et al.} \cite{Liu1992} and has the property that $\mathcal{W}_{JJ}^{\text{vol,L}}=\kappa/J$. This property is useful in that it leads to a constant instantaneous bulk modulus, which in turn provides a constant bulk wave speed (see Equation (\ref{eqn:bulk_wave_speed})). The first derivative of $\mathcal{W}^{\text{vol,q}}$ is linear in $J,$ whereas the first derivative of $\mathcal{W}^{\text{vol,L}}$, and most other formulations for the volumetric energy, is nonlinear in $J$. A more comprehensive list of possible forms for $\mathcal{W}^{\text{vol}}$ is provided by Rossi \textit{et al.} in \cite{Rossi2016}.

The shear modulus, $\mu$, and the bulk modulus, $\kappa$, are proportionality constants relating shear and volumetric energy to their corresponding strains, respectively. They are both derived from the material properties of the solid being modeled. In this paper, $\mu$ and $\kappa$ are defined using Young's modulus, $E$, and Poisson's ratio, $\nu$, via \cite{Lvov2022}
\begin{equation}\label{Moduli}
    \mu = \frac{E}{2(1+\mu)}\qquad \qquad \kappa = \frac{E}{3(1-2\nu)}.
\end{equation}
It is clear from this formulation that $\kappa\rightarrow\infty$ as $\nu\rightarrow0.5$.

\subsection{Dynamics}\label{Dynamics}
This manuscript focuses on numerically solving the Lagrangian conservation of linear momentum in a continuous domain, $\Omega_0$
\begin{equation}
\begin{aligned} \label{CoLM}
    \rho_0(\textbf{X})\dot{\mathbf{V}}(\textbf{X},t)&=\nabla_0\cdot\mathbb{P} (\textbf{X},t)+\rho_0(\textbf{X})\mathbf{B}(\textbf{X},t) & & &\mathbf{X}&\in\Omega_0\\
    \mathbf{\bm{\chi}}(\textbf{X},0)&=\mathbf{\bm{\chi}}_0(\mathbf{X})& & &\mathbf{X}&\in\Omega_0\\
    {\textbf{V}}(\textbf{X},0)&=\textbf{V}_0(\mathbf{X})& & &\mathbf{X}&\in\Omega_0\\
    \mathbf{\bm{\chi}}(\textbf{X},t)&=\textbf{g}_D(\mathbf{X},t)& &&\mathbf{X}&\in\partial\Omega_0^D\\
    \mathbb{P}(\textbf{X},t)\cdot\mathbf{N},(\mathbf{X},t)&=\mathbf{T}(\mathbf{X},t)& &&\mathbf{X}&\in\partial\Omega_0^N
\end{aligned}
\end{equation}
in which $\rho_0$ is the density of the material in reference configuration and $\mathbf{B}$ is a vector of body forces. Dot notation $(\dot{\quad})$ indicates a time derivative and $\nabla_0\cdot(\cdot)$ represents the divergence operator in reference coordinates. $\partial\Omega_D$ and $\partial\Omega_N$ respectively indicate the portions of the boundary of $\Omega$ where Dirichlet and natural boundary conditions are enforced. For simplicity, we assume the deformations on $\Omega_0$ are both isotropic and adiabatic.

By substituting Equations (\ref{Wdev}) and (\ref{Wvol}) into Equation (\ref{PK1}), Equation (\ref{CoLM}) becomes 
\begin{equation}\label{fullCOLM}
    \rho_0  \dot{\mathbf{V}} =  \mu\nabla_0\cdot \left[ J^{-2/3}\left(\mathbb{F} - \frac{1}{3}(\mathbb{F}:\mathbb{F})\mathbb{F}^{-T}\right) \right] + \nabla_0\cdot\left[ \mathcal{W}^{\text{vol}}_J\mathbb{H}\right] +\rho_0\mathbf{B},
\end{equation}

The main challenge in numerically solving Equation (\ref{fullCOLM}) arises from the value of $\kappa$, which approaches infinity as $\nu \rightarrow 0.5$ (see Equation (\ref{Moduli})). This indicates that dynamics in hyperelastic materials become increasingly numerically stiff in the incompressible limit \cite{Leveque2007}, which becomes a necessary consideration when choosing a time integration scheme for solving Equation (\ref{CoLM}). Deformations in $\Omega$ will exhibit both shear waves, which travel at a speed
$$c_\nu = \sqrt{\frac{\mu}{\rho_0}},$$ and bulk waves, which travel at a speed 
\begin{equation}\label{eqn:bulk_wave_speed}
    c_\kappa = \sqrt{\frac{\tilde{\kappa} + 4/3\mu}{\rho_0}},
\end{equation}
in which $\tilde{\kappa}\equiv\kappa J\mathcal{W}^{\text{vol}}_{JJ}$ represents the instantaneous bulk modulus. By the Courant-Friedrichs-Lewy condition \cite{Courant1967}, these wave speeds are crucial for setting constraints on time-step sizes for conditionally stable time integration schemes \cite{Rossi2016,Gil2014}. Because employing the Liu model for the volumetric energy produces a constant instantaneous bulk modulus, this work will use it implicitly when calculating wave speeds and consequent time step size restrictions. $J\approx1$ for incompressible simulations, so these approximations will hold for the quadratic volumetric energy.

\section{Time integration scheme}\label{timestepping}

For nearly- or fully-incompressible problems, the upper limit on explicit time step size is dictated by the speed of bulk waves in the system, the general strategy behind all of these methods is to introduce a Lagrange multiplier pressure field that enforces incompressibility. This pressure field can be updated implicitly, alleviating the time step size restriction associated with bulk waves. Equation (\ref{CoLM}) becomes the system
\begin{align} \label{CoLMsystem}
    \rho_0 \dot{\mathbf{V}} &= \nabla_0\cdot \mathbb{P}^{\text{dev}}(\mathbf{U}) + \nabla_0\cdot (p\mathbb{H}) +\rho_0\mathbf{B}\\
    p &= \frac{\mathcal{W}^{\text{vol}}}{\partial J}(\mathbf{U}) =\mathcal{W}^{\text{vol}}_J.\label{eqn:pressure}
\end{align}
 For the rest of the work, both sides of Equation (\ref{eqn:pressure}) will be divided by $\kappa$. This serves as a way to eliminate the numerical stiffness inherent to the volumetric strain energy; instead of the right hand side of the pressure equation going to infinity, the left hand side goes to zero, reducing to the standard incompressibility condition.

\subsection{Weak formulation and spatial discretization}\label{SpatialDiscretization}
Equations (\ref{CoLMsystem}) and (\ref{eqn:pressure}) can be expressed in weak form by multiplying through by the test functions $\textbf{w}$ and $q$, where $$\textbf{w}\in \mathcal{V}\subseteq H^1(\bar\Omega)^d\qquad q\in \mathcal{P}\subseteq L^2(\bar{\Omega}),$$
and integrating by parts, providing
\begin{equation}\label{WeakForm}
    \begin{aligned}
        (\mathbf{w},\dot{\mathbf{V}})_{\Omega_0} &= (\mathbf{w},\mathbf{T})_{\partial\Omega_0}-(\nabla_0\mathbf{w},\mathbb{P})_{\Omega_0}+ (\mathbf{w},\mathbf{B})_{\Omega_0}\\
        \left(q,\frac{1}{\kappa}p\right)_{\Omega_0}&=\left(q,\frac{1}{\kappa}\mathcal{W}^{\text{vol}}_{J}\right)_{\Omega_0}\\
    \end{aligned}
\end{equation}
in which $(\cdot,\cdot)_{\Omega_0}$ indicates the standard $L^2$ inner product across the undeformed domain. By assuming that the test function space and the trial function space are the same, solving the weak form of Equations (\ref{CoLMsystem}) and (\ref{eqn:pressure}) can be restated as:
\begin{itquote}
    Find $\mathbf{U},\mathbf{V}\in \mathcal{V}$ and  $p\in\mathcal{P}$ such that Equation (\ref{WeakForm}) is satisfied for all $\textbf{w}\in\mathcal{V}$ and $q\in\mathcal{P}$. 
\end{itquote} 

The domain $\Omega$ can be discretized into a triangulation $\mathcal{T}$ composed of elements $K$ and nodes $\textbf{x}_\mathcal{T}$,  such that $\Omega = \bigcup_{K\in \mathcal{T}} K$. $\mathcal{T}$ is equipped with finite element spaces $\mathcal{V}_h$ and $\mathcal{P}_h$ such that the shape functions $\phi_\textbf{v}\in\mathcal{V}_h$ and $\phi_p\in\mathcal{P}_h$ interpolate $\textbf{U},$ $\textbf{V},$ and $p$ at the \textit{N} nodes of $\mathcal{T}$. Approximate solution vectors can be defined such that 
\begin{equation*}    \mathbf{U}(\mathbf{X})\approx\sum^N_{i} \phi_{\textbf{U},i}(\textbf{X})\{\textbf{U}_h\}_i\qquad 
    \mathbf{V}(\mathbf{X})\approx\sum^N_{i} \phi_{\textbf{V},i}(\textbf{X})\{\textbf{V}_h\}_i\qquad 
    p(\mathbf{X})\approx\sum^N_{i} \phi_{p,i}(\textbf{X})\{p_h\}_i.
\end{equation*}

\subsection{Semi-implicit backwards differentiation formula}

The schemes featured in this work implement a modified form of the second-order semi-implicit backward differentiation formula (SBDF2), as presented by Chow \textit{et al.} \cite{Chow2021}. It is often used to solve ordinary differential equations of the form
\begin{equation*}
    \dot{y} = f(y)+g(y),
\end{equation*}
in which $f$ and $g$ respectively represent non-stiff and stiff, potentially-nonlinear operators and $y$ is a function that solves the ODE. The SBDF2 scheme is defined as:
\begin{equation}
    \frac{3y^{n+1} - 4 y^{n} + y^{n-1}}{2\Delta t} = \tilde{f}(y^{n,n-1}) + g(y^{n+1}) = 2f(y^n) - f(y^{n-1}) + g(y^{n+1}), \label{SBDF2}
\end{equation}
Tildes ($\tilde{\cdot}$) denote a second-order Adams-Bashforth (AB2) forward extrapolation in time \cite{Ascher1995, Bashforth1883}. SBDF2 is capable of damping high frequency error terms while maintaining stability if the eigenvalues corresponding to $g$ are sufficiently larger than those associated with $f$ \cite{Chow2021,Ascher1995, Xu2006}. It has been utilized in the context of finite element methods but not specifically to problems in incompressible elasticity. Decomposing Equations (\ref{CoLMsystem}) and (\ref{eqn:pressure}) into a set of first order PDEs and applying Equation (\ref{SBDF2}) yields
 \begin{equation}
     \begin{aligned}
         \rho_0\frac{3\textbf{V}^{n+1} - 4 \textbf{V}^{n} + \textbf{V}^{n-1}}{2\Delta t} &= \nabla_0^h\cdot\left(\tilde{\mathbb{P}}^{\text{dev}}(\mathbf{U}^{n,n-1})\right) + \nabla_0^h\cdot\left(p^{n+1}\mathbb{H}(\mathbf{U}^{n+1})\right) + \rho_0 \mathbf{B}^{n+1}\\
         \frac{3\textbf{U}^{n+1} - 4 \textbf{U}^{n} + \textbf{U}^{n-1}}{2\Delta t} &= \mathbf{V}^{n+1}\\
         \frac{1}{\kappa}p^{n+1} &=\frac{1}{\kappa}\mathcal{W}^{\text{vol}}_{J}(\mathbf{U}^{n+1}),
     \end{aligned}
 \end{equation}
in which $\textbf{V}=\dot{\textbf{U}}$ is the velocity field defined on $\Omega_0$. In this formulation, updating the velocity field depends on an implicit evaluation of both the pressure and the cofactor, which is a nonlinear function of $\textbf{U}$. One of the goals of the present research is to avoid nonlinear solvers and their associated computational cost, so having an implicit evaluation that depends on both displacement and pressure in the same term is suboptimal.

\subsection{Modified semi-implicit backward differentiation formula}\label{MSBDF2section}
The goal of implementing a semi-implicit method, such as SBDF2, is to minimize the amount of numerical complexity required to solve a given PDE. With regard to nearly- and fully- incompressible elasticity, the scheme should be designed to solve for the pressure field implicitly while treating every other relevant term explicitly. Consequently, the cofactor in the volumetric term of the PK1 stress tensor can also be extrapolated using an additional AB2 approximation, providing the following discretized system:
\begin{equation}\label{MSBDF2}
     \begin{aligned}
         \rho_0\frac{3\textbf{V}^{n+1} - 4 \textbf{V}^{n} + \textbf{V}^{n-1}}{2\Delta t} &= \nabla_0^h\cdot\left(\tilde{\mathbb{P}}^{\text{dev}}(\mathbf{U}^{n,n-1})\right) + \nabla_0^h\cdot\left(p^{n+1}\tilde{\mathbb{H}}(\mathbf{U}^{n,n-1})\right) + \rho_0 \mathbf{B}^{n+1}\\
         \frac{3\textbf{U}^{n+1} - 4 \textbf{U}^{n} + \textbf{U}^{n-1}}{2\Delta t} &= \mathbf{V}^{n+1}\\
        \frac{1}{\kappa} p^{n+1} &= \frac{1}{\kappa} \mathcal{W}^{\text{vol}}_{J}\left(J(\mathbf{U}^{n+1})\right).
     \end{aligned}
 \end{equation}
 Strictly speaking, utilizing an explicit extrapolation in both terms of the PK1 stress tensor makes this method distinct from SBDF2, and consequently is referred to as Modified SBDF2, or MSBDF2. In contrast to SBDF2, MSBDF2 is most useful being applied to problems of the form 
 $$ \dot{y}(t)=f(y,t)+h(y,t)g(y,t)$$ where $f$ and $g$ are both non-stiff operators, and $h$ is a stiff operator.
    
The incompressibility condition is enforced through the solution to the pressure equation, which may be nonlinear in $J$ depending on the chosen form of $\mathcal{W}^{\text{vol}}$. Setting up the fully discretized scheme in the following section is dependent on linearizing both $\mathcal{W}^{\text{vol}}$ in $J$ and $J$ in $\textbf{U}$. The incompressiblity condition is thus weakly enforced for large deformations if the linearization of $\mathcal{W}^{\text{vol}}$ becomes inaccurate. 

\subsection{Forward Euler / backward differentiation formula }
While SBDF2 is a well studied integrator, its application to hyperbolic problems is not well studied, and consequently there is no guarantee that MSBDF2 is going to be a suitable algorithm for solving PDEs like Equation (\ref{CoLM}). Adams-Bashforth extrapolations feature harsher time step size restrictions than other methods \cite{Bashforth1883}. Consequently, using another method of extrapolation, like forward Euler, may introduce superior stability properties. 

A forward Euler / Backward Differentiation Formula Runge-Kutta scheme can be expressed as 
\begin{equation}\label{FEBDF2}
     \begin{aligned}
     \mathbf{U}^*&=\mathbf{U}^n+\Delta t \mathbf{V}^n \\
         \rho_0\frac{3\textbf{V}^{n+1} - 4 \textbf{V}^{n} + \textbf{V}^{n-1}}{2\Delta t} &= \nabla_0^h\cdot\left({\mathbb{P}}^{\text{dev}}(\mathbf{U}^*)\right) + \nabla_0^h\cdot\left(p^{n+1}{\mathbb{H}}(\mathbf{U}^*)\right) + \rho_0 \mathbf{B}^{n+1}\\
         \frac{3\textbf{U}^{n+1} - 4 \textbf{U}^{n} + \textbf{U}^{n-1}}{2\Delta t} &= \mathbf{V}^{n+1}\\
        \frac{1}{\kappa} p^{n+1} &= \frac{1}{\kappa} \mathcal{W}^{\text{vol}}_{J}\left(J(\mathbf{U}^{n+1})\right),
     \end{aligned}
 \end{equation}
where the first stage is denoted by an asterisk. This algorithm will be referred to as FEBDF2 for the remainder of the manuscript. The difference in computational cost between the two algorithms is negligible.
  \subsection{Fully discretized scheme:}
 By combining the spatial discretization provided in section (\ref{SpatialDiscretization}) and  MSBDF2 or FEBDF2, Equation (\ref{CoLMsystem}) can be discretely represented in the following block linear system:
\begin{equation}\label{LinearSystem}\left[
\begin{array}{cc}
 \textbf{M}_\textbf{V}    & \textbf{M}_{\textbf{V}p} \\
 \textbf{M}_{p\textbf{V}}    & \textbf{M}_{p}
\end{array} \right]   \left[\begin{array}{c}
     \mathbf{V}_h^{n+1} \\
     p_h^{n+1}
\end{array}\right] = \left[\begin{array}{c}
    \mathbf{R}_\textbf{V} \\
    \textbf{R}_p
\end{array}\right],
\end{equation}
with mass matrix contributions of
\begin{align*}
\mathbf{M}_{\textbf{V},ij} &= \frac{\rho_0}{\Delta t}\left(\phi_{\textbf{V},i}(\textbf{X}),\phi_{\textbf{V},j}(\textbf{X})\right)_{\Omega_0}& \mathbf{M}_{\textbf{V}p,ij} = \frac{2}{3}\left(\nabla_0^h \phi_{\textbf{V},i}(\textbf{X}),\bar{\mathbb{H}}\phi_{p,j}(\textbf{X})\right)&_{\Omega_0}\\
    \mathbf{M}_{p\textbf{V},ij} &= -\frac{\Delta t}{\kappa} \left(\phi_{p,i}(\textbf{X}),\mathcal{W}_{JJ}^{vol,n}\mathbb{H}^n:\nabla_0^h\phi_{\textbf{V},j}(\textbf{X})\right)_{\Omega_0}&
    \mathbf{M}_{p,ij} = \frac{1}{\kappa}\left(\phi_{p,i}(\textbf{X}),\phi_{p,j}(\textbf{X})\right)&_{\Omega_0}
\end{align*}
and residual contributions of 
\begin{align*}
    \mathbf{R}_{\textbf{V},i}&=\frac{\rho_0}{3\Delta t}\left(\phi_{\textbf{V},i}(\textbf{X}),4\textbf{V}_h^n-\textbf{V}_h^{n-1}\right)_{\Omega_0}+\frac{2}{3}(\phi_{\textbf{V},i}(\textbf{X}),\mathbf{T}^{n+1})_{\partial\Omega_0}-\frac{2}{3}(\nabla_0^h\phi_{\textbf{V},i}(\textbf{X}),\bar{\mathbb{P}}^{\text{dev}})_{\Omega_0}\\
    \mathbf{R}_{p,i}&=\frac{1}{\kappa}\left(\phi_{p,i}(\textbf{X}),\mathcal{W}^{\text{vol}}_J\left(J(\textbf{U}_h^{n})\right)\right)_{\Omega_0}
\end{align*}
Bars $(\bar{\cdot})$ are used to denote either extrapolations via AB2 extrapolations for MSBDF2 or forward Euler extrapolations for FEBDF2. It is imperative to note that $||\textbf{M}_p||\rightarrow0$ for $\kappa\rightarrow\infty$, which establishes this system as a saddle-point problem in the incompressible limit \cite{Brezzi1974}. Following Kadapa \cite{Kadapa2021} and standard techniques for solving the incompressible Stokes equations \cite{Wang2020}, the solution to this system can be updated by solving a Schur complement system,
\begin{equation}\label{SchurComplement}
    S = (\mathbf{M}_p-\mathbf{M}_{p\textbf{V}}\mathbf{M}^{-1}_{\textbf{V}}\mathbf{M}_{\textbf{V}p}),
\end{equation}
via: 
\begin{equation}\label{Iteration}
    \begin{aligned}
       1) \quad p_h^{n+1} &= S^{-1}\left(R_p-\mathbf{M}_{p\textbf{V}}\mathbf{M}_{\textbf{V}}^{-1}R_\textbf{V}\right)\\
       2)\quad \textbf{V}_h^{n+1}&=\mathbf{M}_{\textbf{V}}^{-1}\left(R_\textbf{V}-\mathbf{M}_{\textbf{V}p}p_h^{n+1}\right)\\
       3) \quad \textbf{U}_h^{n+1} &= 1/3\left(4\textbf{U}_h^n-\textbf{U}_h^{n-1}+2\Delta t \textbf{V}_h^{n+1}\right)
    \end{aligned}
\end{equation}

The inversion of the Schur complement is analogous to solving a weighted Helmholtz problem, and in the limit $\kappa\rightarrow\infty,$ it is analogous to solving a weighted Poisson problem.  

 \section{Analysis of MSBDF2 and FEBDF2}\label{Analysis}
\subsection{Local truncation error analysis}
As established by Taylor analysis in Appendix \ref{appendix_A}, the local truncation error for MSBDF2 at time $t^n$ can be expressed as
\begin{equation}
    \mathcal{L}^n = \Delta t^3\left[\frac{\rho_0}{12}\dddot{\textbf{V}}^n +\nabla_0\cdot\left(\ddot{\mathbb{P}}^{\text{dev},n} +\frac{1}{2}\ddot{p}^n\mathbb{H}^n-\dot{p}^n\dot{\mathbb{H}}^n-\frac{1}{2}p^n\ddot{\mathbb{H}}^n\right)\right] + O(\Delta t^4).
\end{equation}
This establishes MSBDF2 as a second order method in time. Similarly, the local truncation error for FEBDF2 can be expressed as 
\begin{equation}
    \mathcal{L}^n=\Delta t^3\left[\frac{\rho_0}{12}\dddot{\textbf{V}}^n -\nabla_0\cdot\left(\frac{1}{2} \mathbf{V}^n\cdot\partial_{\mathbf{UU}}{\mathbb{P}}^{\text{dev},n}\cdot\mathbf{V}^n +\frac{1}{2}\ddot{p}^n\mathbb{H}^n+\dot{p}^n\dot{\mathbb{H}}^n+\frac{1}{2}p^n\mathbf{V}^n\cdot\partial_{\mathbf{UU}}{\mathbb{H}}^n\cdot\mathbf{V}^n\right)\right] + O(\Delta t^4)
\end{equation}
establishing FEBDF2 as a second order method in time. The order of accuracy of both schemes is heavily dependent on the extrapolation of $\mathbb{H}$; without this extrapolation, both methods reduce to first order accuracy.

\subsection{Stability} \label{sec:stability}
Numerical integration schemes for nonlinear PDEs often lack concrete stability properties, unlike their linear counterparts. A common way to study the stability of a nonlinear numerical method is by analyzing its effect on the total energy of the system. A stable method should satisfy the condition 
\begin{equation}\label{EnergyCondtion}
    E^{n+1}\leq E^{n},
\end{equation}
where $E$ represents the total energy present within the computational domain \cite{Xu2019,Simo1992}.

The equation for the conservation of linear momentum is hyperbolic by the nature of time-dependent conservation laws, and thus should preserve a constant energy given a lack of external forcing \cite{Tadmor2012}. A natural energy for the system is provided by the sum of the kinetic energy and the internal energy of the system as defined by the strain energy density function:
\begin{equation}
   0\leq E(t) = \int_{\Omega_0} \left(\frac{\rho_0}{2}(\mathbf{V}\cdot\mathbf{V}) + \mathcal{W}^{\text{dev}}(\mathbb{F})+\frac{1}{2}p(J-1)\right)dV.    \label{Energy}
\end{equation}
The volumetric term of the strain energy density function was expressed using the quadratic formulation and a substitution of the definition of pressure from Equation (\ref{CoLMsystem}) allows for the energy to be approximated in the fully incompressible case. Boundary forces are neglected in this representation, based on the assumption of no external forcing; computing energy in a system with outside forcing would necessitate the inclusion of boundary forces. 

\subsubsection{Linear Stability: MSBDF2}

The goal of stability analysis in general is to compute a maximum time step size, $\Delta t_{max}$, such that the energy of the system does not grow without external forcing \cite{Leveque2007}. For problems in the nonlinear regime, establishing energy stability is very difficult for integration schemes lacking a particular structure \cite{Simo1992,Kang2023}. Methods for establishing formal energy stability for BDF-based integration schemes applied to hyperbolic PDEs do not appear to have been developed. The following analysis circumvents the nonlinearity of $\mathcal{W}$ by studying incompressible deformations in a linear regime. The following time step analysis does not guarantee stability for all nonlinear problems, but it does provide a suitable upper bound on time step size for a large class of problems. 

Because Equation (\ref{LinearSystem}) is a decomposed hyperbolic PDE, and thus features simultaneous updates to velocity, displacement, and pressure, it is sensible to represent the evolution of the system via the evolution of a system of difference equations. Consequently, the following stability analysis focuses on the eigenvalues of the matrix that defines the evolution of the discrete system.

When all nonlinear operators are replaced by their linearized forms, Equation (\ref{LinearSystem}) takes the form 
\begin{equation}
\begin{aligned}   \rho_0\mathbf{M}_{\textbf{v}} \mathbf{V}^{n+1} &= 4/3\textbf{M}_{\textbf{v}}  \mathbf{V}^n - 1/3\mathbf{M}_{\textbf{v}}  \mathbf{V}^{n-1} - 2/3 \Delta t\left( \textbf{K}^{\text{dev}}(2\mathbf{U}^{n}-\textbf{U}^{n-1})+\textbf{K}^{\text{vol}}(\textbf{U}^{n+1})\right)\\
   \mathbf{M}_{\textbf{v}} \textbf{U}^{n+1} &= 4/3\mathbf{M}_{\textbf{v}} \textbf{U}^{n} - 1/3 \mathbf{M}_{\textbf{v}} \textbf{U}^{n-1} + 2/3 \Delta t \textbf{M}_{\textbf{v}} \mathbf{V}^{n+1},
\end{aligned}
\end{equation}
   where $\textbf{K}^{\text{dev}}$ and $\textbf{K}^{\text{vol}}$ are the linear operators representing the weak form of the deviatoric and the volumetric components to the PK1 stress tensor respectively. Restructuring the system such that it can be written in terms of consecutive iterations of a solution vector,
   $$\mathbf{U}_0^{n+1} = \textbf{U}^{n}_1 \qquad \quad \mathbf{V}_0^{n+1} = \textbf{V}^{n}_1,$$
yields 
\begin{equation}\label{difference_system}
\mathbf{A}_1\mathbf{Y}^{n+1}=\textbf{A}_0\textbf{Y}^n,\end{equation}
in which

$$
\begin{array}{lr}
    \mathbf{A}_1 =\left[
\begin{array}{cccc}
 \textbf{M}_\textbf{v}    & 0 & 2/3\Delta t\textbf{K}^{\text{vol}} & 0  \\
   0  & \textbf{M}_\textbf{v} & 0 & 0 \\
    -2/3\Delta t\textbf{M}_\textbf{v}  & 0 & \textbf{M}_\textbf{v} & 0 \\
     0  & 0 & 0 & \textbf{M}_\textbf{v}
\end{array} \right]   & \textbf{Y}^{n+1}=\left[\begin{array}{c}
     \mathbf{V}_1^{n+1} \\
     \mathbf{V}_0^{n+1} \\
     \mathbf{U}_1^{n+1} \\
     \mathbf{U}_0^{n+1}
\end{array}\right] \\ 
\textbf{A}_0=  \left[
\begin{array}{cccc}
 4/3\textbf{M}_\textbf{v}    & -1/3\textbf{M}_\textbf{v}  & -4/3\textbf{K}^{\text{dev}} & 2/3\Delta t \textbf{K}^{\text{dev}}  \\
   \textbf{M}_\textbf{v}  & 0 & 0 & 0 \\
    0  & 0 & 4/3\textbf{M}_\textbf{v} & -1/3\textbf{M}_\textbf{v} \\
     0  & 0 & \textbf{M}_\textbf{v} & 0
\end{array} \right]  &
\textbf{Y}^{n}=\left[\begin{array}{c}
     \mathbf{V}_1^{n} \\
     \mathbf{V}_0^{n} \\
     \mathbf{U}_1^{n} \\
     \mathbf{U}_0^{n}
\end{array}\right]
\end{array}
$$

This system can be simplified by approximating the effect of $K^{\text{dev}}$ and $K^{\text{vol}}$ via eigenvalues. 
Following Kadapa and Joly \cite{Kadapa2021, Joly2007}, the approximation $$||\textbf{M}_\textbf{v}^{-1}\textbf{K}^{\text{dev}}||\leq\frac{4\mu}{\rho_0h^2}\equiv \lambda$$ can be used to approximate the maximum deformation due to deviatoric forces, where $h$ is introduced as the minimum length between vertices of an element. There is no directly analogous approximation for volumetric deformations, so the approximation
$$0<||\mathbf{M}_\textbf{v}^{-1}\textbf{K}^{\text{vol}}||\leq c$$
is used, where $c$ is left unknown. The effect of solving Equation (\ref{difference_system}) can be  determined by studying
\begin{equation}
\big\|\textbf{Y}^{n+1}\big\|\leq\big\|A_1^{-1}A_0\big\|\textbf{ }\big\|\textbf{Y}^n\big\|,
\end{equation}
where $\big\|A_1^{-1}A_0\big\|$ can be approximated by
$$ \max_\omega\left(\left[\begin{array}{cccc}
 1    & 0 & 2/3\Delta tc & 0  \\
   0  & 1 & 0 & 0 \\
    -2/3\Delta t  & 0 & 1 & 0 \\
     0  & 0 & 0 & 1
\end{array} \right]^{-1}\left[
\begin{array}{cccc}
 4/3    & -2/3  & -4/3\Delta t\lambda & 2/3\Delta t \lambda  \\
   1  & 0 & 0 & 0 \\
    0  & 0 & 4/3 & -1/3 \\
     0  & 0 & 1 & 0
\end{array} \right]\right).$$
Solving the quartic characteristic polynomial attributed to it and solving for the $\Delta t$ that guarantee $\max(\omega)\leq 1$ may be technically possible, but infeasible. Consequently, given a value for $\lambda$, a maximum time step size can be numerically approximated such that all eigenvalues of the system are bounded by 1.

\begin{figure}[t!]
    \centering    \includegraphics[width=0.99\linewidth]{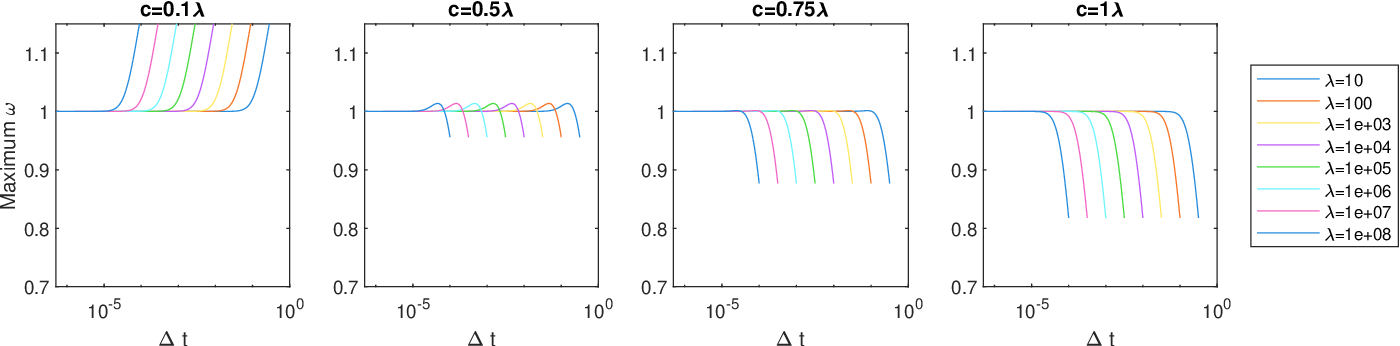}
    \caption{Maximum $\omega$ obtained for varying values of $c$, $\lambda$, and $\Delta t$. Values for $c$ were approximated as scalar multiples of $\lambda$. }
    \label{fig:MSBDF2_dt_eig}
\end{figure}

This style of analysis highlighted a fascinating result; MSBDF2's stability is completely dependent on the behavior of the volumetric component of the energy, denoted by $c$. As displayed in Figure \ref{fig:MSBDF2_dt_eig}, if $c\geq\sim0.75\lambda,$ $|\omega|<1$, and stability is ensured for all time step sizes. If $c<\sim0.75\lambda,$ $\omega$ may be greater than 1, but converges to 1 if $\Delta t$ is small enough. Therefore, if the effect of enforcing the volumetric constraint is proportionally large, the method will be unconditionally stable, but if the effect of enforcing the incompressibility constraint is trivial, then the method is unstable with the caveat that a small $\Delta t$ will minimize the effect of this instability. 

It is impossible to approximate a maximum time step size that would guarantee stability for a particular problem; the value for $c$ evolves with the numerical solution, which implicitly makes it dependent on the time step size. Consequently, a time step size restriction for MSBDF2 can only be approximated experimentally. 

\begin{figure}[t!]
    \centering
    \includegraphics[width=0.6\linewidth]{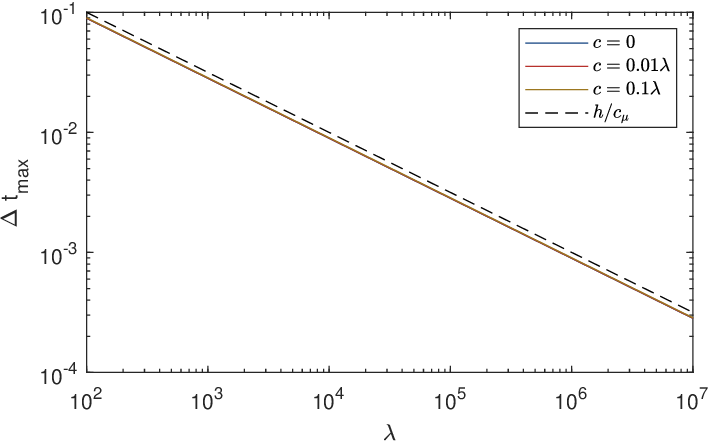}
    \caption{Maximum time step size, $\Delta t_{max}$ for FEBDF2 calculated with various values for $\lambda$ and $c$.}
    \label{fig:FEBDF2_max_dt}
\end{figure}

\subsubsection{Linear Stability: FEBDF2}
The same analysis outlined above extends to FEBDF2 as well. Equation (\ref{FEBDF2}) can be expressed in a linear form via 
\begin{equation}
\begin{aligned}   \rho_0\mathbf{M}_{\textbf{v}} \mathbf{V}^{n+1} &= 4/3\textbf{M}_{\textbf{v}}  \mathbf{V}^n - 1/3\mathbf{M}_{\textbf{v}}  \mathbf{V}^{n-1} - 2/3 \Delta t\left( \textbf{K}^{\text{dev}}(\mathbf{U}^{n}+\Delta t\textbf{V}^{n-1})+\textbf{K}^{\text{vol}}(\textbf{U}^{n+1})\right)\\
   \mathbf{M}_{\textbf{v}} \textbf{U}^{n+1} &= 4/3\mathbf{M}_{\textbf{v}} \textbf{U}^{n} - 1/3 \mathbf{M}_{\textbf{v}} \textbf{U}^{n-1} + 2/3 \Delta t M_{\textbf{v}} \mathbf{V}^{n+1}.
\end{aligned}
\end{equation}
After employing the same substitutions used in the previous section, stability for the system is guaranteed by ensuring 

$$ \max_\omega\left(\left[\begin{array}{cccc}
 1    & 0 & 2/3\Delta tc & 0  \\
   0  & 1 & 0 & 0 \\
    -2/3\Delta t  & 0 & 1 & 0 \\
     0  & 0 & 0 & 1
\end{array} \right]^{-1}\left[
\begin{array}{cccc}
 2/3(2-\Delta t^2 \lambda)    & -1/3  & -2/3\Delta t\lambda & 0  \\
   1  & 0 & 0 & 0 \\
    0  & 0 & 4/3 & -1/3 \\
     0  & 0 & 1 & 0
\end{array} \right]\right)$$
is less than 1, where $\omega$ are the eigenvalues of the system. Though similar in form to MSBDF2's stability condition, this system has radically different behavior. As shown in Figure \ref{fig:FEBDF2_max_dt}, the maximum time step size that guarantees stability for FEBDF2 scales like $h/c_{\mu},$ which is the expected time step size restriction for a semi-implicit method \cite{Kadapa2021, Gil2014,Lahiri2005}. For simulations in this work, FEBDF2 was run with a CFL number of 0.5, which was sufficient to guarantee stability for nonlinear problems \cite{Courant1967}.

\section{Computational setup}

The efficiency of MSBDF2 and FEBDF2 is partially dictated by the inversion of the $\textbf{M}_\mathbf{V}$ matrix. The finite element basis functions used in this work partition unity; this property can be taken advantage of by defining a new mass matrix with diagonal values defined as the sum of all the values along their respective row in $\textbf{M}_\mathbf{V}$:
\begin{equation}\label{Mass matrix lumping}
    \mathbf{M}^{LM}_{\mathbf{V},ii}=\rho_0\sum_j\int_{\Omega_0}\phi_{\mathbf{V},i}(\mathbf{X})\phi_{\mathbf{V},j}(\mathbf{X})dV.
\end{equation}
The code implemented in this work employs mass lumping to trivialize the inversion of the mass matrix in computing the Schur complement via Equation (\ref{SchurComplement}), significantly reducing the computational cost of each iteration of the entire scheme. The error produced by diagonalization is not higher than error introduced through spatial discretization, and consequently does not negatively affect the accuracy of the overall method\cite{Hughes2000,Voet2023}.

As discussed earlier, Equation (\ref{LinearSystem}) reduces to a saddle-point system in the incompressible limit, and consequently Taylor-Hood\cite{Taylor1973} element pairs are utilized for the joint velocity / pressure finite element spaces to satisfy the Ladyzhenskaya-Babuska-Brezzi condition\cite{Brezzi1974}. For this work, quadratic elements were used for the velocity space, and linear elements were used for the pressure space. For quadrilateral and hexahedral elements, this corresponds to the Q2/Q1 finite element pairing. In contrast, mass lumping with simplices equipped with a quadratic approximation space leads to non-positive values along the diagonal of the mass matrix. This corresponds to the solid being modeled having non-positive mass, which is non-physical. This is remedied by enriching the simplex finite-element space approximating the velocity space with a bubble function \cite{Brezzi1992,Cohen2001}. Introducing a single degree of freedom within the center of the element (and in the center of each face in 3D) enables the computation of a positive-definite lumped mass matrix. Quadratic finite element spaces equipped with an additional centroid degree of freedom on simplices are referred to as $P2+$, so the enriched Taylor-Hood finite element for simplex elements is $P2+/P1$.   

The conjugate gradient method\cite{Hestenes1952} with a diaonal preconditioner is used for the inversion of any consistent mass matrix if necessary, and the lumped velocity mass matrix does not need a preconditioner due to the triviality of inverting a diagonal matrix. The most computationally expensive aspect of MSBDF2 is the inversion of the Schur complement operator. Due to the asymmetry of the operator, the generalized minimum-residual method is used for the Schur complement solve. The Schur complement is a discretized version of a weighted Helmholtz/Poisson operator, so an algebraic multigrid preconditioner \cite{Falgout2006} can be initialized with a preconditioning matrix, $P$, that is the discretized weak form of the weighted Helmholtz/Poisson operator:
\begin{equation}\label{Preconditioner}
    P^n_{ij}=\int_{\Omega_0}\frac{1}{\kappa}\phi_{p,i}\phi_{p,j} -\frac{2\Delta t }{3} \tilde{\mathbb{H}}^{n,n-1}\nabla_0\phi_{p,i}\cdot\mathbb{H}^n\nabla\phi_{p,j}dV.
\end{equation}
Assembling this preconditioning matrix incurs minimal cost because it reuses terms already computed during system matrix assembly. It also changes form depending on the value of $\kappa$, so it is a suitable preconditioner for both nearly- and fully- incompressible problems.

\section{Numerical Tests}
Numerical studies reported herein used a code written using the \textit{deal.II} finite element library \cite{dealii}, with regular meshes composed of $Q2/Q1$ element pairs. Results were visualized using ParaView \cite{paraview} and convergence plots were constructed using $MATLAB$ \cite{matlab}.

Because many of the examples studied in this section do not have analytical solutions in Lagrangian coordinates, the standard definition for error for an arbitrary solution vector $y$, 
\begin{equation}\label{error}
    \epsilon_y^n = y^n-y(t^n),
\end{equation}
cannot be used for establishing convergence rates. Convergence rates are instead determined by studying the difference between numerical solutions obtained before and after refining the time step size, i.e.
\begin{equation}
    \varepsilon_{y,2\Delta t}^N=y^N_{2\Delta t}-y^N_{\Delta t},
\end{equation}
where $N$ is used to denote the solution was obtained at the final time step. For instances where an analytical solution is available, the notation "$\epsilon$" is used, and if not, "$\varepsilon$" is used. 

The utilized norms were left unspecified in the errors defined above. In this work, three norms are used to measure accuracy / convergence: 
\begin{equation}
    \begin{aligned}
        \big\|\varepsilon\big\|_{L^1}&=\int_\Omega |\varepsilon|dV \\
        \big\|\varepsilon\big\|_{L^2}&=\left(\int_\Omega |\varepsilon|^2 dV\right)^{1/2} \\
        \big\|\varepsilon\big\|_{L^\infty}&= \inf\left\{C:C\geq|\varepsilon|\text{ for almost all $\textbf{X}\in\Omega$}\right\}.
    \end{aligned}
\end{equation}
The behavior of these different norms can vary. For example, the existence of a persistent singularity in the solution can cause a lack of convergence in the $L^\infty$ norm, but the other two norms will converge if the singularity is located at a single point. It is consequently beneficial to use all three to properly capture convergence. 

\begin{figure}[b!]
    \centering
    \begin{tikzpicture}
    \node (A){{\includegraphics[width=0.9\linewidth]{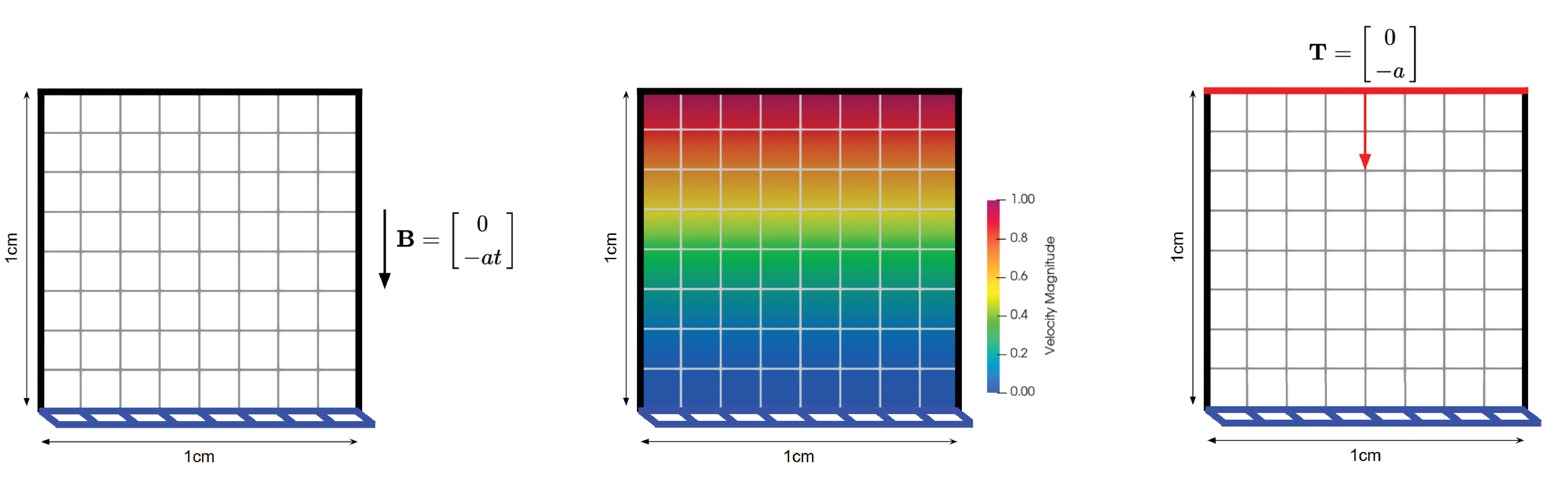}}};
    \path (A.south west) -- (A.south east) node[pos=.14] {\footnotesize(a)} node[pos=.51] {\footnotesize(b)} node[pos=.86] {\footnotesize(c)};
    \end{tikzpicture}
    \caption{Force diagrams for the three unit square simulations. Blue edges indicate homogeneous Dirichlet boundary conditions. (a) depicts a unit square equipped with a constant body force with magnitude $a=25\text{g cm}/
     \text{s}^3$. (b) depicts a  equipped with an initial velocity field defined by Equation (\ref{initial_velocity}). (c) depicts a unit square with a downward constant traction applied to the top face, where $a=-2.5\text{dyn}/\text{cm}^2$.}
    \label{fig:square_diagram}
\end{figure}
\begin{table}[t!]
    \centering
    \begin{tabular}{|c| c c c|}
    \hline
        $\nu$ & 0.4 & 0.49 & 0.5\\
        \hline
       $\mu$ $(dyn/cm^2)$ & 35.714 & 33.557 & $33.\bar{3}$ \\ 
       $\kappa$ $(dyn/cm^3)$ & 3.3$\times10^2$ &3.3$\times10^3$ & 3.3$\times10^4$ \\
       $c_\mu$ $(cm/s)$ & 5.79 & 5.793 & 5.774 \\
       $c_\kappa$ $(cm/s)$ & 19.518 & 58.12& $\infty$ \\
        \hline
        $\Delta t_{\mu}$ $(s)$& 0.0105 & 0.0108 & 0.0108 \\
        $\Delta t_{\kappa}$ $(s)$& 0.0032 & 0.0011 & 0\\
        \hline
    \end{tabular}
    \caption{Table of material parameters obtained from $E=100
    \text{dyn}/\text{cm}^2$, $\rho_0=1g/cm^3$, and various values for $\nu$, as well as the minimum time step sizes require to resolve shear waves ($\Delta t_\mu=h_{min}/c_\mu$) and bulk waves ($\Delta t_\kappa=h_{min}/c_\kappa$). }
    \label{tab:US_table}
\end{table}
\subsection{Deformed unit square} \label{sec:US}

To validate the predicted rate of convergence, a 2D square with dimensions $1\text{cm}\times 1\text{cm}$, $E=100\text{dyn/cm}^2$, and density $1\text{g/cm}^3$ was placed under 3 different forms of deformation field: a uniform body force, an initial velocity field, and a constant applied traction. Error convergence due to temporal refinement is calculated for varying values of $\nu$.

\begin{figure}[b!]
\centering
\begin{tikzpicture}
\node(A){\includegraphics[width=0.9\textwidth]{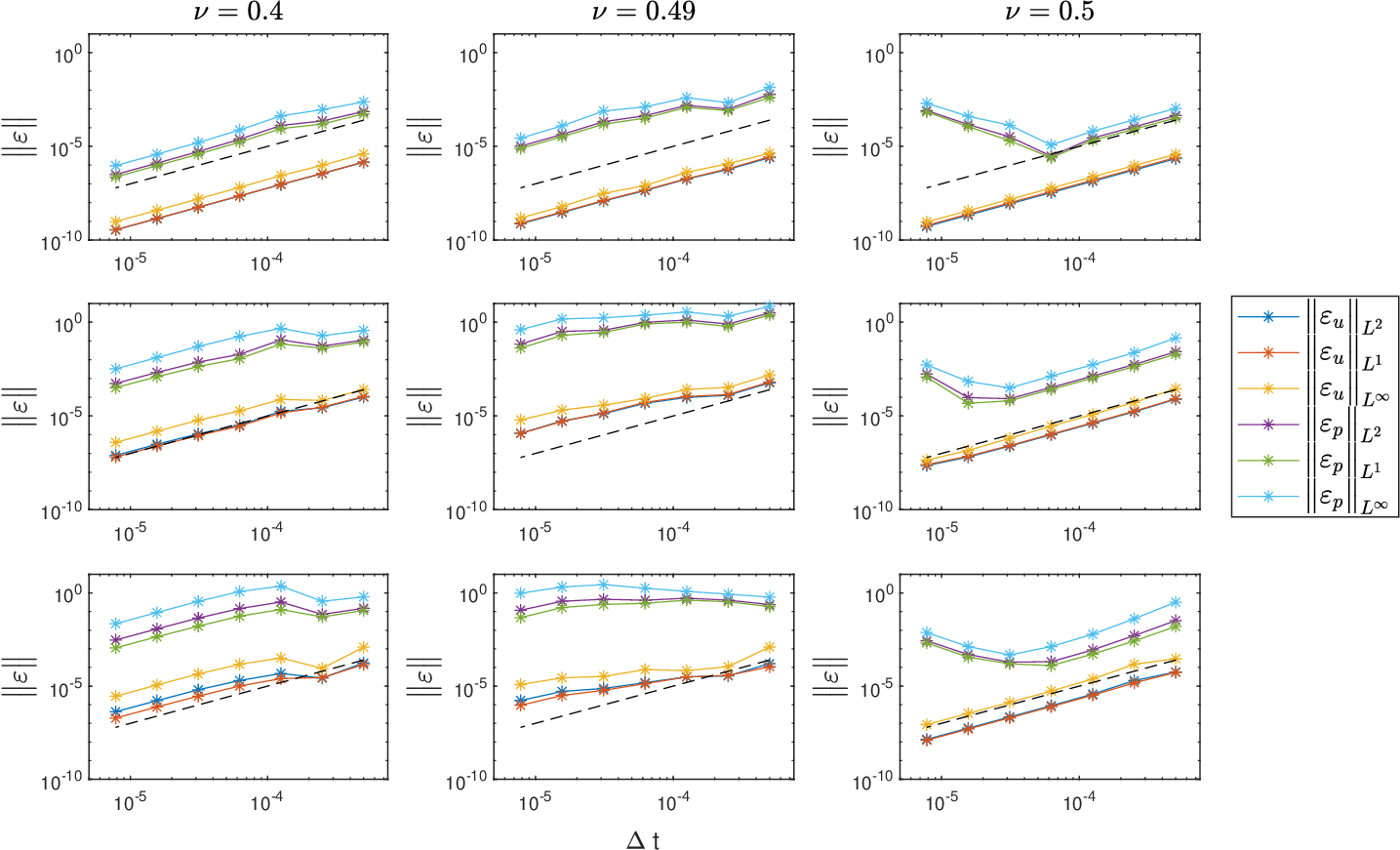}
};
\path (A.south west) -- (A.north west) node[pos=.22, above, rotate=90]
{Tr}
node[pos=.53, above, rotate=90] {IV}
node[pos=.84, above, rotate=90] {BF};
\end{tikzpicture}
\caption{Convergence in behavior for deformations of a unit square as computed by MSBDF2. Each row depicts $\big\|\varepsilon\big\|_{L^p}$ values obtained from body force, initial velocity, and traction respectively.}\label{fig:MSBDF2_convergence}
\end{figure}

\subsubsection{Body force}

The first numerical test is a unit square under a uniform body force that increases linearly in time (see Figure \ref{fig:square_diagram}). The bottom face of all meshes are equipped with homogeneous Dirichlet boundary conditions, and all other faces are traction-free unless specific. The body force is given by $\mathbf{B}=[0,-at]$. The force coefficient, $a$, is set to $25\text{g cm/s}^3$ for these simulations. 

The first row of Figure (\ref{fig:MSBDF2_convergence}) demonstrates that MSBDF2 is capable of obtaining second order convergence for $\nu=0.4$ and $\nu=0.49,$ though a loss of convergence emerges in the pressure field for $\nu=0.5$. The source of this lack of convergence is unknown, though it could derive from the emergence of artificial numerical shocks. 

\begin{figure}[t!]
\centering
\begin{tikzpicture}
\node(A){\includegraphics[width=0.9\textwidth]{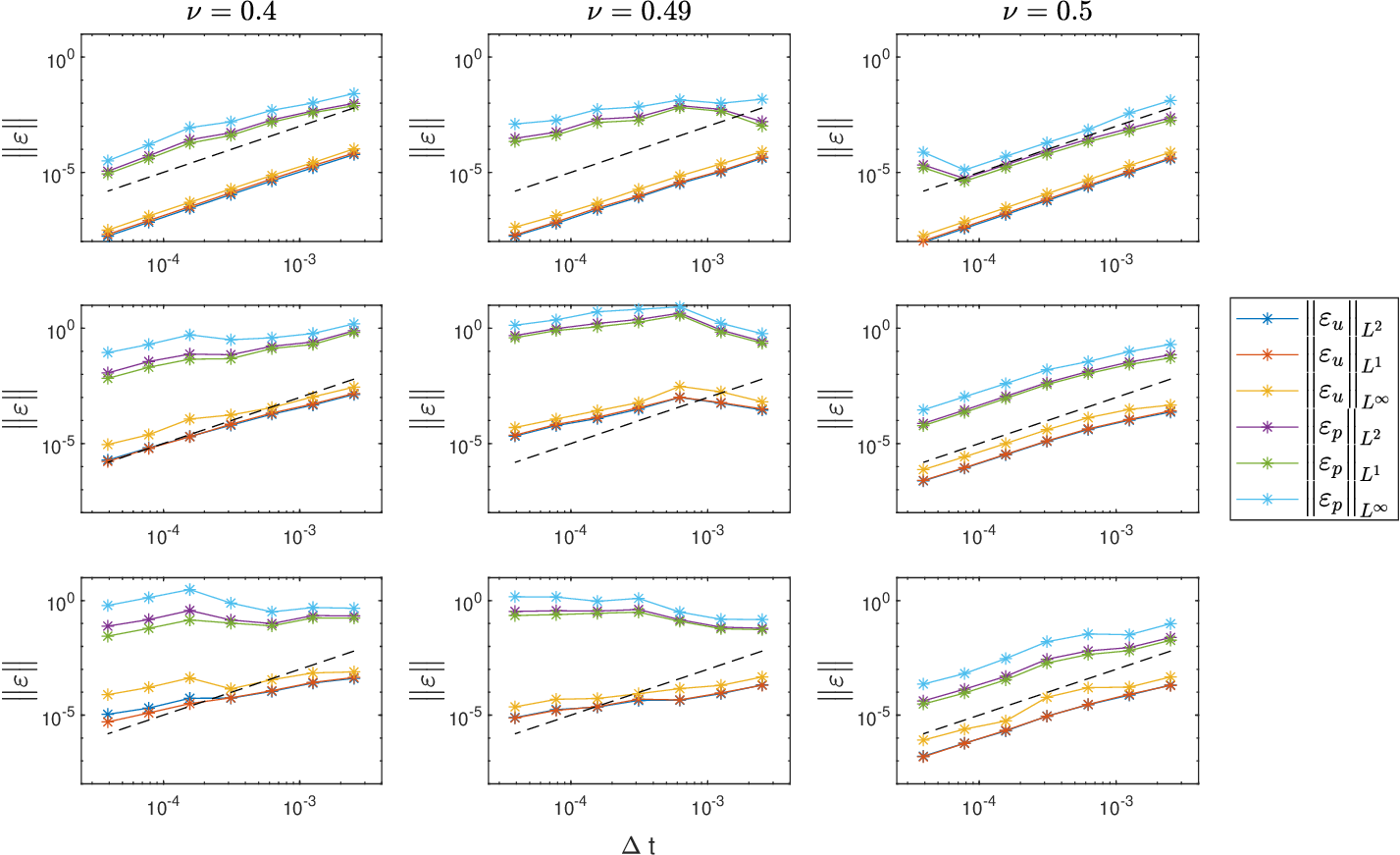}
};
\path (A.south west) -- (A.north west) node[pos=.22, above, rotate=90]
{Tr}
node[pos=.53, above, rotate=90] {IV}
node[pos=.84, above, rotate=90] {BF};
\end{tikzpicture}
\caption{Convergence in behavior for deformations of a unit square as computed by FEBDF2. Each row depicts $\big\|\varepsilon\big\|_{L^p}$ values obtained from body force, initial velocity, and traction respectively.}\label{fig:FEBDF2_convergence}
\end{figure}

The first row of Figure (\ref{fig:FEBDF2_convergence}) demonstrates that FEBDF2 is capable of obtaining second order convergence for $\nu=0.4$ and $\nu=0.5,$ though a similar lack of proper convergence is observed in the pressure field in the $\nu=0.5$ case. Sub-quadratic convergence is observed in the pressure field when $\nu=0.49$.

\subsubsection{Initial velocity field}
The second numerical test studies the convergence obtained when the unit square is equipped with the initial velocity field $\label{initial_velocity}
    \mathbf{V}_0=\left[0,-\sin(\pi/2 X_2)\right].$
The magnitude of the initial velocity field is depicted in Figure  (\ref{fig:square_diagram}). The mesh has traction-free boundary conditions on all faces barring the bottom face, which has homogeneous Dirichlet boundary conditions. This initial velocity field is not divergence-free, indicating that there should be an initial shock when run in the nearly- and fully- incompressible regime.

The second row of Figure (\ref{fig:MSBDF2_convergence}) demonstrates that MSBDF2 is capable of obtaining second order convergence for $\nu=0.4$ and $\nu=0.49,$ though a loss of convergence emerges in the pressure field when $\nu=0.5$.

The second row of Figure (\ref{fig:FEBDF2_convergence}) demonstrates that FEBDF2 is capable of obtaining second order convergence for $\nu=0.4$ and $\nu=0.5$. Subquadratic convergence is observed in the pressure field when $\nu=0.49$, though only for adequately small time step sizes.

\subsubsection{Constant traction}

The second test run is a unit square with a traction boundary condition applied to the top face. As shown in Figure (\ref{fig:square_diagram}), the faces on the sides of the domain are traction-free, and the bottom face is equipped with homogeneous Dirichlet boundary conditions. The traction vector is defined by 
$\mathbf{T} = \left[
    0 ,
     -2.5\text{dyn/cm}^2
\right].$
The mesh is unloaded and undeformed at $t=0,$ so the traction vector induces a volumetric shock that, given the nature of the governing PDE, should propagate at $c_\kappa$. 

The third row of Figure (\ref{fig:MSBDF2_convergence}) demonstrates that MSBDF2 is capable of obtaining second order convergence for $\nu=0.4$ and in the velocity field for $\nu=0.5$ though a loss of convergence emerges in the pressure field when $\nu=0.5$. Sub-quadratic convergence is observed in both the displacement and pressure fields for $\nu=0.49.$ 

The third row of Figure (\ref{fig:FEBDF2_convergence}) demonstrates that FEBDF2 produces sub-quadratic convergence for $\nu=0.4$ and in the displacement field for $\nu=0.49,$ though a lack of convergence is observed in the pressure field for $\nu=0.49.$ Second order accuracy is obtained in the fully incompressible case. 

\begin{table}[t!]
    \centering
    \begin{tabular}{|c| c c c|}
    \hline
        $\nu$ & 0.4 & 0.49 & 0.5\\
        \hline
       $\mu$ $(dyn/cm^2)$ & 892.86 & 838.93 & $833.\bar{3}$ \\ 
       $\kappa$ $(dyn/cm^3)$ & $8.\bar{3}e3$ &$8.\bar{3}e4$ & $\infty$ \\
       $c_\mu$ $(cm/s)$ & 94.49 & 91.59 & 91.29 \\
       $c_\kappa$ $(cm/s)$ & 97.59 & 290.6& $\infty$ \\
       \hline
    \end{tabular}
    \caption{Table of material parameters obtained from $E=2500\text{dyn/cm}^2$, $\rho_0=0.1\text{g/cm}^3$, and various values for $\nu$.}
    \label{tab:cook_table}
\end{table}
\begin{figure}[b!]
    \centering
    \subfloat[]{
    \includegraphics[width=0.4\linewidth]{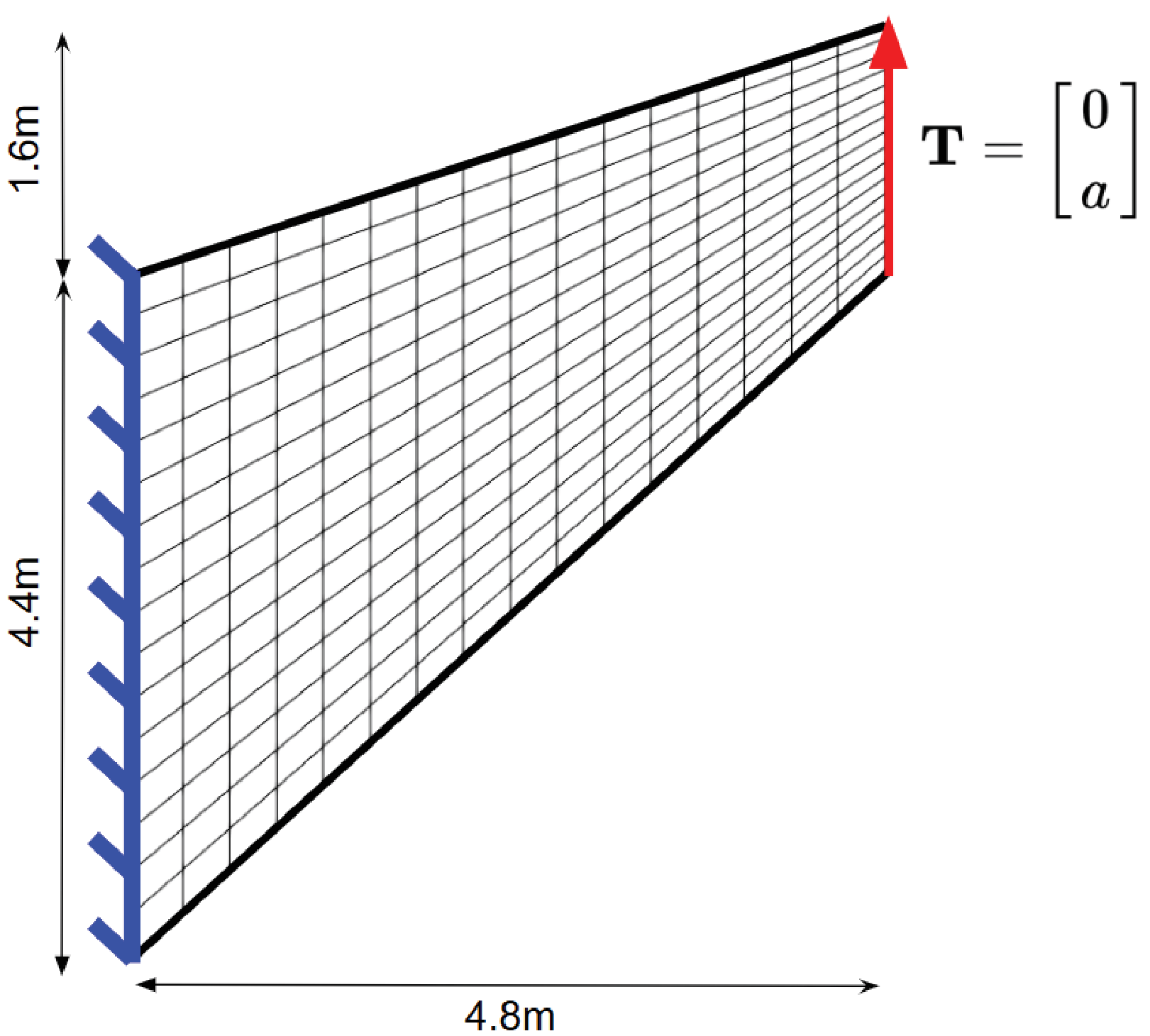}}
    \subfloat[]{
    \includegraphics[width=0.28\linewidth]{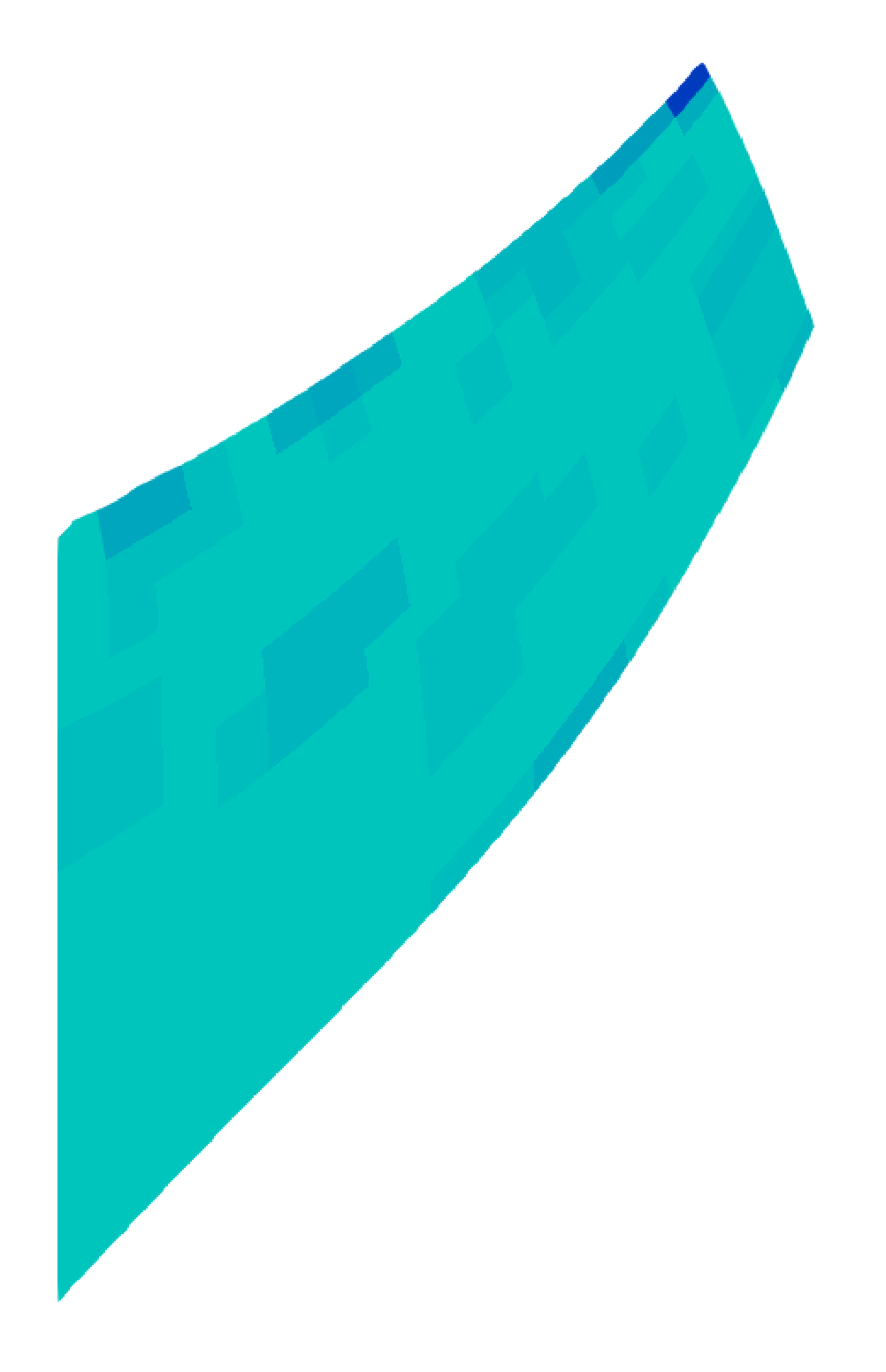}}
    \subfloat[]{
    \includegraphics[width=0.28\linewidth]{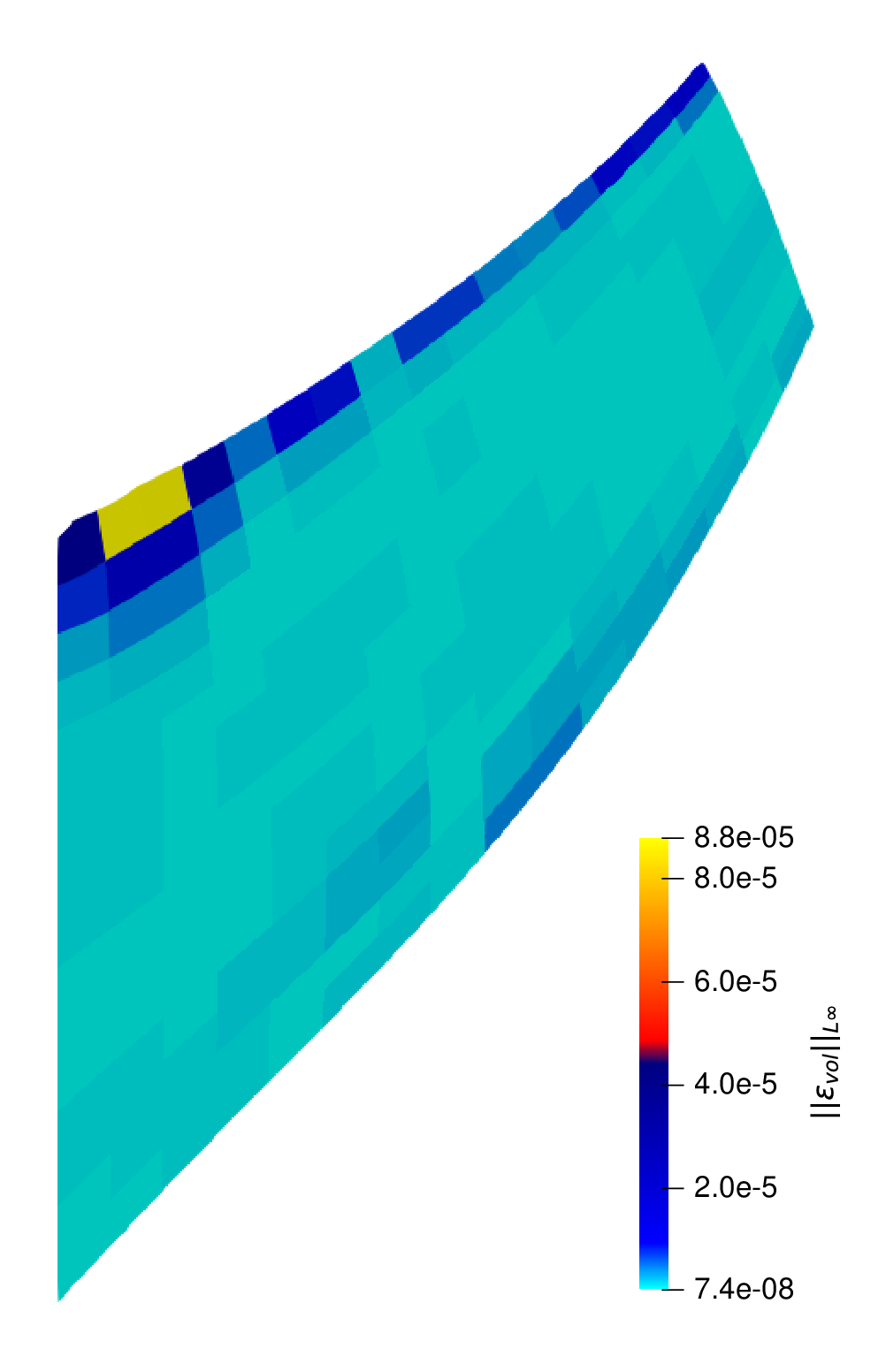}}
    \caption{(a) A force diagram for the Cook's Membrane problem. Black edges represent traction-free boundaries, blue edges indicate homogeneous Dirichlet boundary conditions, and red edges indicate faces with enforced traction. For the simulations in this work, $a=62.5\text{dyn/cm}^2$. Subfigures (b) and (c) show maximum volumetric error per cell after 15s of simulation time using $\mathcal{W}^{\text{vol,L}}$ and $\mathcal{W}^{\text{vol,L}}$ respectively.}
    \label{fig:Cooks_membrane}
\end{figure}

\subsection{Cook's membrane}\label{CooksMembrane}
The second numerical test used to validate MSBDF2 is Cook's membrane, a classic benchmark problem for algorithms in computational solid mechanics \cite{Rossi2016, Kadapa2021}. The problem involves a trapezoidal mesh with an upward-pointing traction applied to the right face (see Figure \ref{fig:Cooks_membrane}a). The membrane has a Young's modulus of $E=2500\text{dyn/cm}^3$, a density of $0.1\text{g/cm}^3$, and various values for the Poisson ratio. The material properties derived from these parameter values are summarized in Table \ref{tab:cook_table}.

This is a standard benchmark because the prescribed traction causes the formation of a pressure singularity at the upper left corner. Simulation results for this benchmark are displayed in Figure \ref{fig:Cooks_membrane}. The functionality that this specific experiment sought to study was MSBDF2's ability to conserve volume. A model for an incompressible solid should naturally preserve volume throughout the computational domain. It is impossible for volume to be conserved pointwise, but the total volume, given by
$ V(t) = \int_{\Omega_t} dv,$
should remain constant.

\begin{figure}[t!]
    \centering
    \begin{tikzpicture}
    \node (A) {{\includegraphics[width=0.9\linewidth]{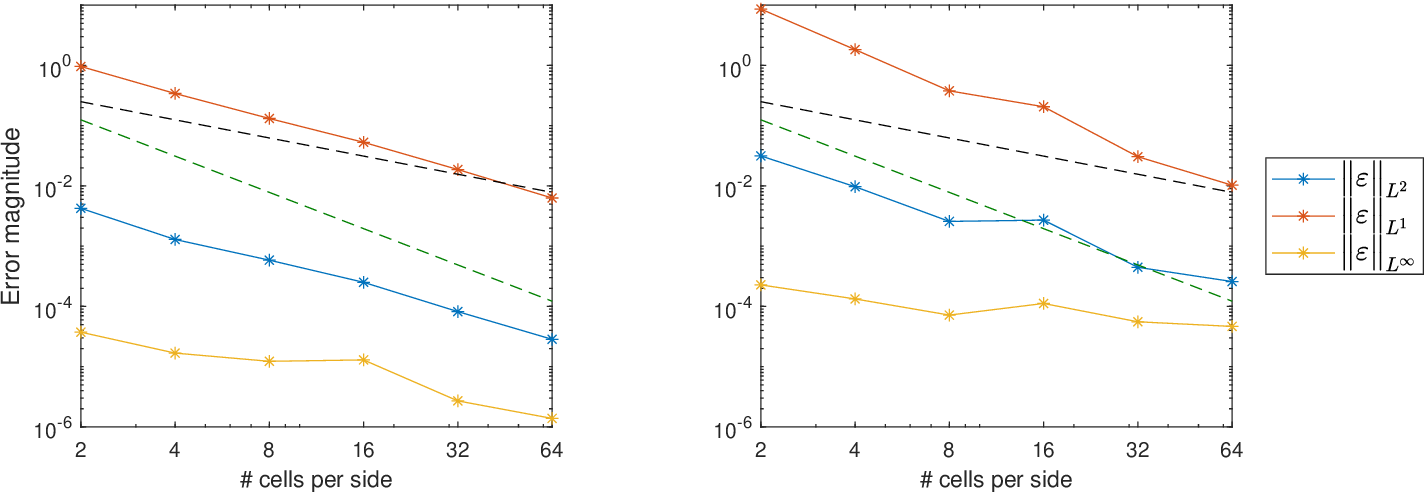}}};
    \path (A.south west) -- (A.south east) node[pos=.23] {\footnotesize(a)} node[pos=.71] {\footnotesize(b)} ;
    \end{tikzpicture}
    \caption{Convergence plots for the volumetric error in Cook's membrane simulations run with FEBDF2 using $\mathcal{W}^{\text{vol,q}}$ in Subfigure (a) and $\mathcal{W}^{\text{vol,L}}$ in Subfigure (b). A linear trendline is provided in black and a quadratic trendline is provided in green.}
    \label{fig:cook_volume_error}
\end{figure}

\begin{figure}[b!]
    \centering
    \includegraphics[width=0.9\linewidth]{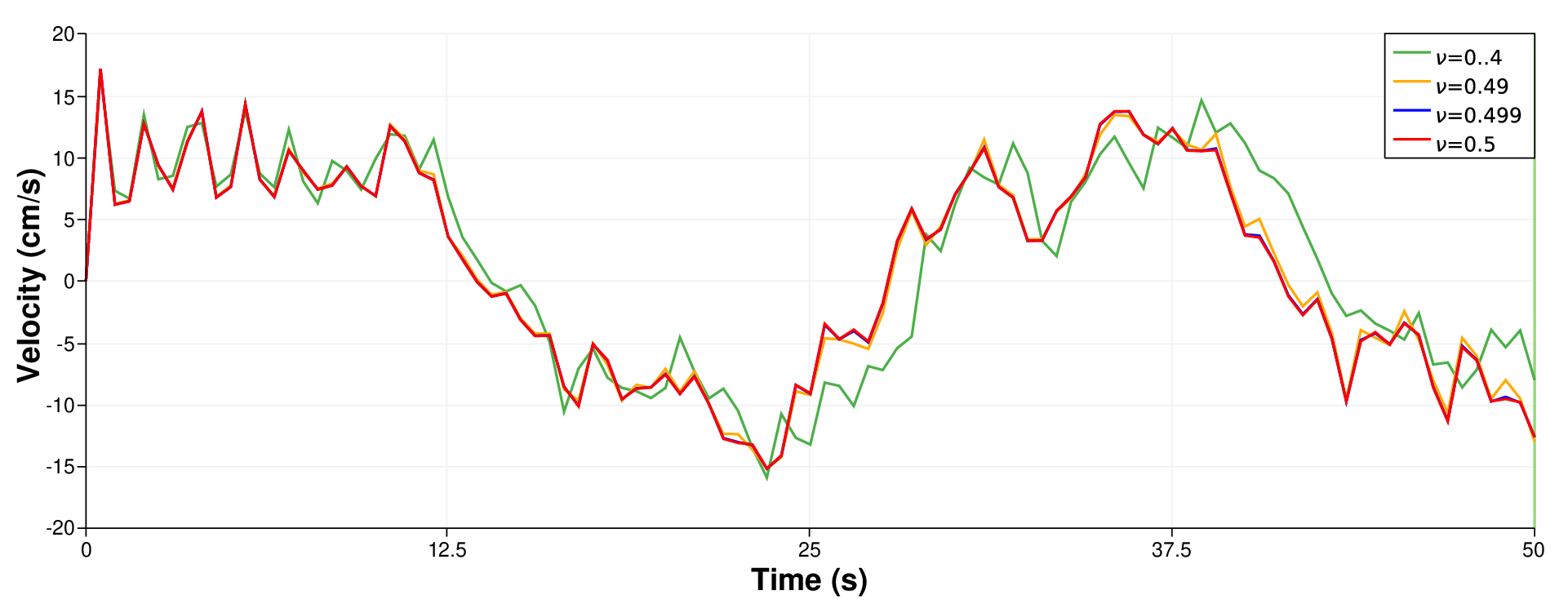}
    \caption{Velocity at the top right tip of the cook's membrane, as calculated by FEBDF2. Data displayed is temporally under-resolved due to rapid oscillations obfuscating the convergence in behavior as $\nu\rightarrow0.5$.}
    \label{fig:tip_velocity}
\end{figure}

It is possible to derive an error that can be used to validate the volume conservation properties of the integration scheme by invoking the properties of $J$. Because the Jacobian of deformation is defined throughout the domain, the volumetric error defined as
\begin{equation}
    \epsilon_{\text{vol}}(\textbf{X},t) = \big\|J(\mathbf{X},t)-1 \big\|_{L^p}
\end{equation}
can be evaluated at each time step for little additional cost. The three norms discussed before are again utilized here to identify the presence of a singularity. If the volume of the simulation is indeed preserved, this error should converge to zero under grid refinement.

The main purpose of this numerical test is to highlight whether or not the choice of the formulation for $\mathcal{W}^{\text{vol}}$ affects the ability for MSBDF2 to conserve volume. Figure \ref{fig:cook_volume_error} depicts the effect of simultaneous refinement in both time and space: each refinement entails dividing the grid size by two and the time step size by four. Because volumetric error should be determined by spatial resolution, this imbalance in refinement is intended to minimize the effect of temporal error in the results.

\subsection{Twisting Column}\label{TwistingColumn}
\begin{table}[b!]
    \centering
    \begin{tabular}{|c| c c c|}
    \hline
        $\nu$ & 0.4 & 0.49 & 0.5\\
        \hline
       $\mu$ $(\text{dyn/cm}^3)$ & $4.29\times10^6$ &$4.03\times10^6$ & $4.00\times10^6$ \\
       $\kappa$ $(\text{dyn/cm}^2)$ & $4.00\times10^7$ & $4.00\times10^8$ & $\infty$ \\ 
       $c_\mu$ $(\text{cm/s})$ & 1974 & 1913 & 1907 \\
       $c_\kappa$ $(\text{cm/s})$ & 6761 & 20133& $\infty$ \\
        \hline
        $\Delta t_{\mu}$ $(\text{s})$ & $1.27\times10^{-4}$ & $1.31\times10^{-4}$ & $1.31\times10^{-4}$ \\
        $\Delta t_{\kappa}$ $(\text{s})$&  $3.7\times10^{-5}$ & $1.24\times10^{-5}$ & 0\\
        \hline
    \end{tabular}
    \caption{Table of material and numerical parameters used in the twisting column simulations. These values are obtained from $E=1.2\times10^7\text{dyn/cm}^2$, $\rho_0=1.1\text{g/cm}^3$, various shown values for $\nu$, the characteristic length $h$ is set at $0.5\text{cm}$}
    \label{tab:MSBDF2_column_table}
\end{table}

The results in Figure (\ref{fig:cook_volume_error}) show that the quadratic formulation for $\mathcal{W}^{\text{vol}}$ is more effective at preserving volume under simultaneous temporal and spatial refinement. Volumetric error using the quadratic formulation is superlinear in all norms. The Liu model only provides sub-quadratic convergence in the $L^1$s, while the $L^2$ and $L^\infty$ norms approach zero'th order convergence. This indicates that a singularity emerges in the $J$ field when using the Liu model, which disallows the $\mathcal{W}^{\text{vol,L}}$ from capturing incompressibility as efficiently as the quadratic model under simultaneous spatial and temporal refinement.

\begin{figure}[t!]
    \centering
    \begin{tikzpicture}
    \node (A){{\includegraphics[width=0.85\linewidth]{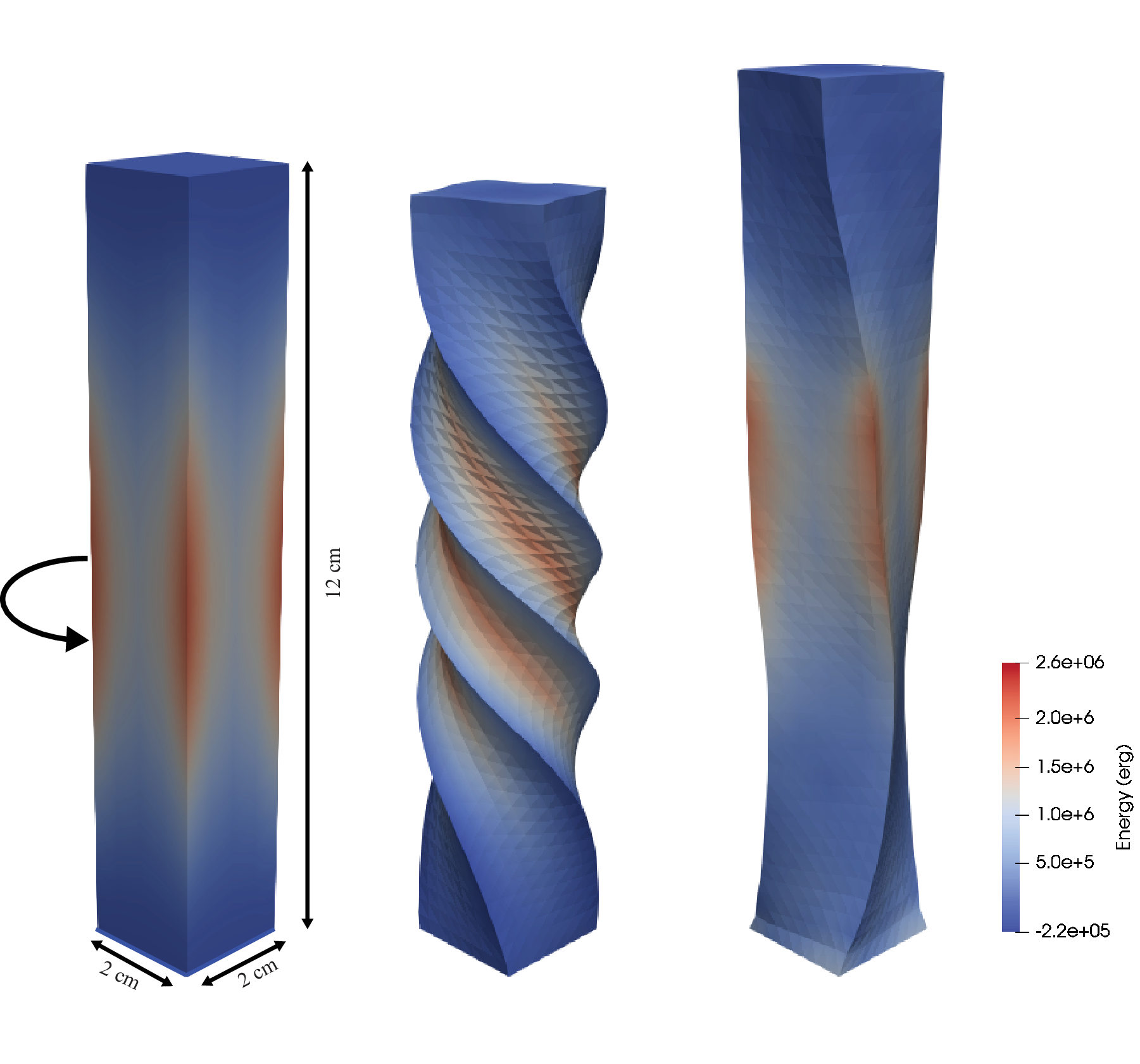}}};
    \path (A.south west) -- (A.south east) node[pos=.17] {\footnotesize(a)} node[pos=.45] {\footnotesize(b)} node[pos=.74] {\footnotesize(c)};
    \end{tikzpicture}
    \caption{(a) The initial state of the twisting column simulation, with the velocity field being defined by Equation (\ref{twist_V0}). (b) and (c) respectively display the result of the fully-incompressible simulation after 0.02 and 0.04 seconds.}
    \label{fig:column_diagram}
\end{figure}

The final numerical test used to validate MSBDF2 is the twisting column. This is also a classical benchmark in solid mechanics because it features very large material parameters and leads to large-scale deformations. Consequently it is a fitting test for the robustness of a numerical scheme. The column is a $2\text{cm}\times2\text{cm}\times12\text{cm}$ rectangular prism with material parameters $E=1.2\times 10^7\text{dyn/cm}^2$ and $\nu =0.4,$ $0.49,$ and $0.5$. The mesh was refined twice, providing element side lengths of $0.25\text{cm}$. The simulation parameters that derive from these parameters are described in Table \ref{tab:MSBDF2_column_table}.

The column is equipped with an initial velocity field described by 
\begin{equation}\label{twist_V0}
    \mathbf{V}(\mathbf{X},t) = \left[\begin{array}{c}
      -1500 \sin\left(\frac{\pi}{12}{X}_3\right) {X_2}\\
      1500 \sin\left(\frac{\pi}{12}{X}_3\right) {X_1} \\
     0 
\end{array}
\right].
\end{equation}
Snapshots of the energy evolution in the column are shown in Figure (\ref{fig:column_diagram}).
\begin{figure}[t]
    \centering
\begin{subfloat}[]{
    \includegraphics[width=0.45\linewidth]{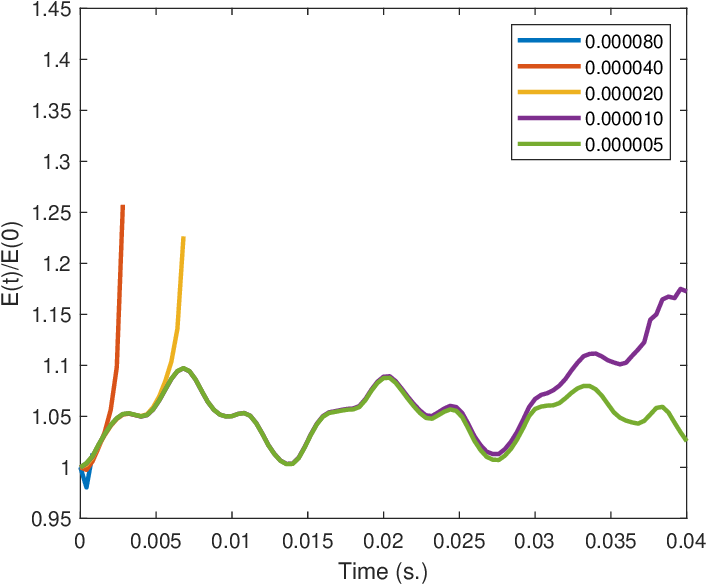}}
\end{subfloat}
\begin{subfloat}[]{
    \includegraphics[width=0.45\linewidth]{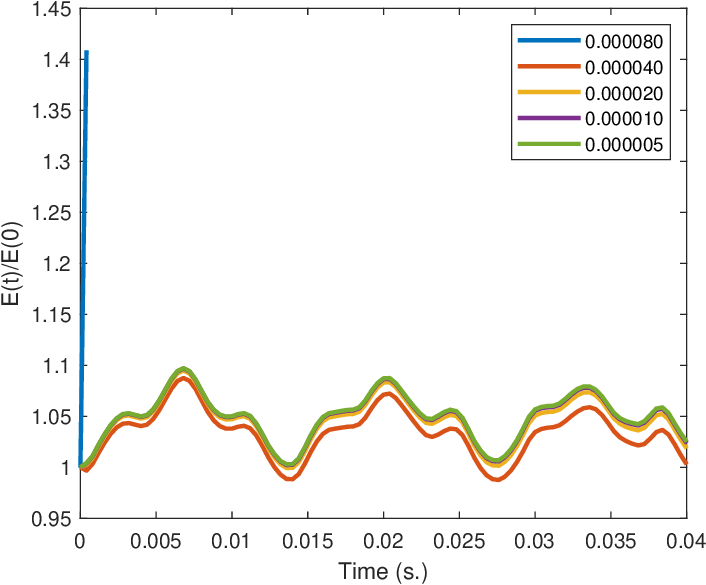}}
\end{subfloat}
    \caption{Evolution of energy in the twisting column simulation as calculated by MSBDF2 (a) and FEBDF2 (b) using various time step sizes. Energy is scaled according to the energy initially present in the system.}
    \label{fig:short_energy}
\end{figure}

The focus of this set of numerical experiments is to demonstrate how MSBDF2 and FEBDF2 preserve energy, defined in Equation (\ref{Energy}). As discussed in Section \ref{sec:stability}, it is ideal to have a notion as to how a numerical method affects the overall energy of a system as it evolves in time. Establishing a proper energy bound on the evolution conservation of linear momentum using the proposed algorithms is not yet possible, though it is possible to measure the evolution of energy through the direct computation of Equation (\ref{Energy}) at each time point. Though the time-step analysis presented earlier in the chapter only applies to linear problems, studying the effect of time step size on energy should provide insight on the stability of MSBDF2/FEBDF2--even in nonlinear problems. 

Figure (\ref{fig:short_energy}) depicts the effect of time step size on the evolution of total energy in a twisting column over time. It is clear that large time step sizes lead to non-physical energy spikes, which in turn cause the simulations to break. Numerical integration schemes for the equations of motion ideally should be able to conserve energy, but the initial conditions for this simulation do not include a pressure field that balances out the prescribe initial kinetic energy, so a bounded energy is expected instead of a non-increasing energy. This expectation is consistent with other results from the literature \cite{Lahiri2005}. Obtaining initial conditions that would demonstrate non-oscillatory energy would require an equilibrium solve or damping, both of which the current state of the MSBDF2/FEBDF2 code is not equipped for.

\begin{figure}[b!]
    \centering
    \includegraphics[width=0.9\linewidth]{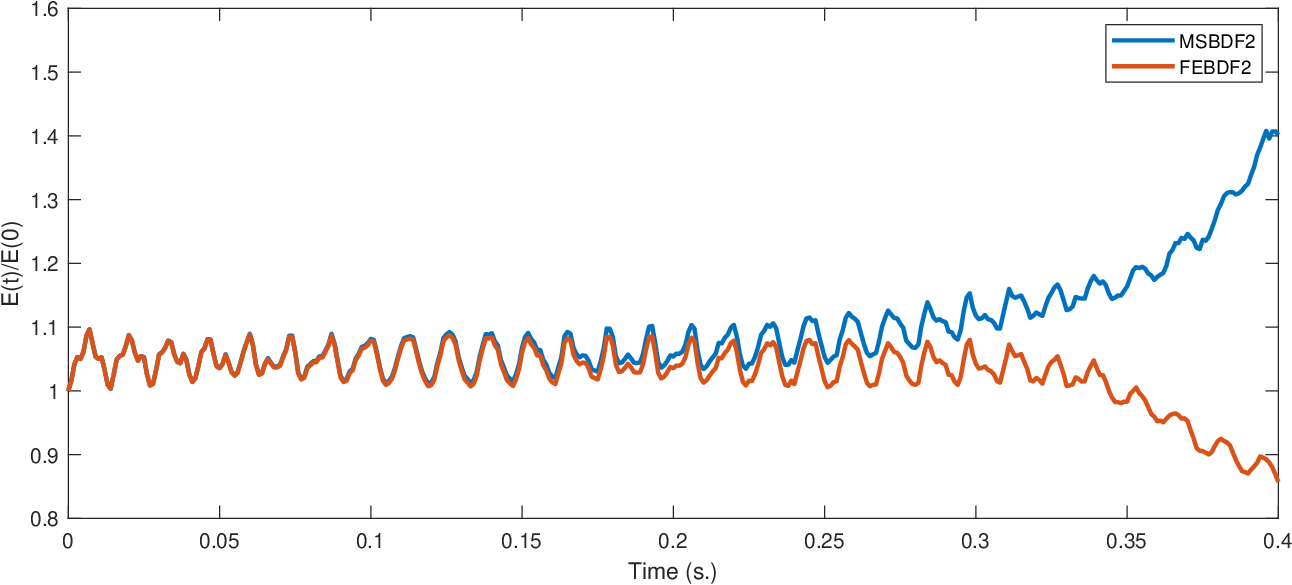}
    \caption{The evolution of energy in an incompressible twisting column run over 0.4 seconds, as calculated by MSBDF2 and FEBDF2. Both simulations were run with $\Delta t =5\times10^{-6}s$}
    \label{fig:long_energy}
\end{figure}

Figure (\ref{fig:short_energy}a) highlights one of the main problems inherent to MSBDF2; not only is it not demonstrably conditionally stable, it is under-dissipative to the point of instability unless extremely small time step sizes are used.  This indicates that even though the maximum eigenvalues in Figure (\ref{fig:MSBDF2_dt_eig}) plateau to 1 as time step size decreases, true convergence may never be achieved regardless of time step size. Figure (\ref{fig:short_energy}b) indicates that FEBDF2 does produces a bounded energy if the time step size less than $\approx0.5\Delta t_{min}$. Convergent solutions with larger time step sizes feature an initial dissipation of energy; this dissipation disappears with smaller time steps.

The long-term evolution of energy in the twisting column simulation can be observed in Figure (\ref{fig:long_energy}). Simultation results were calculated with $\Delta t =5\times10^{-6}\text{s}
$, so both MSBDF2 and FEBDF2 should have produced bounded energies, according to the qualitative results shown in Figure (\ref{fig:short_energy}) and the theoretical results obtained in Section (\ref{sec:stability}). MSBDF2 does not produce a bounded energy, whereas FEBDF2 does, confirming that MSBDF2 can only imitate a notion of stability with small time step sizes, whereas FEBDF2 is a truly stable method. Both simulations develop structural instabilities after enough time, though this type of numerical error is expected.

\section{Discussion and Conclusion}
This paper introduces two integration schemes for resolving large-scale dynamics in nearly- and fully-incompressible solids, both of which are based on the semi-implicit backwards differentiation formula (SBDF2). The introduced schemes modify SBDF2 by multiplicatively decomposing the numerically stiff terms into stiff and nonstiff contributions to reduce computational cost. The MSBDF2 scheme uses an Adams-Bashforth approximation to extrapolate non-stiff terms, whereas the FEBDF2 scheme instead employs a forward Euler stage.

Both of these schemes are proven to be second-order accurate, and the computational results presented in Section (\ref{sec:US}) validate this conclusion, particularly for fully incompressible materials or those with $\nu\approx0.4$. In the nearly incompressible case, results from both MSBDF2 and FEBDF2  converge sub-quadratically if the test problem produces shocks. The convergence of dynamic behavior in Figure (\ref{fig:tip_velocity}) indicates that the lack of convergence in problems involving shocks in nearly-incompressible materials can be circumvented by utilizing a fully-incompressible approximate material model.

Linear stability analysis cannot be used to determine a maximum time-step size that would guarantee that MSBDF2 does not diverge. However, the analytical results implied that MSBDF2 would be able to imitate stability if adequately small time step sizes are used. The same analysis indicated that FEBDF2 is stable if $\Delta t\approx 0.5c_\mu/h$. These predictions were computationally verified in Figures (\ref{fig:short_energy}) and (\ref{fig:long_energy}). These stability results, coupled with second order accuracy, establish FEBDF2 as a very capable integrator for fully incompressible dynamics. MSBDF2, while directly derived from a well-studied semi-implicit integration scheme, does not exhibit similarly favorable stability properties.

While the results in this work establish FEBDF2 as both a stable and accurate integrator, they also highlight a potential weakness inherent to the application of linearized semi-implicit schemes used for problems in nonlinear incompressible elasticity. Methods that depend on the linearization of the incompressibilty constraint necessarily do not perfectly preserve volume. However, as demonstrated in Figure (\ref{fig:cook_volume_error}), utilizing a formulation for the volumetric error such that $\mathcal{W}^{\text{vol}}_{JJ}=\kappa$ ensures greater volume conservation upon spatial and temporal refinement. 

The inability of the Liu formulation to conserve volume is clearly not due to the choice of material model; it is a problem rooted in the derivation of linearized semi-implicit schemes. MSBDF2 and FEBDF2 depend on the linearization of the pressure equation, which has little effect when using $\mathcal{W}^{\text{vol,q}}$ because it is already linear in $J$. The pressure term derived from the Liu model is not linear in $J$, and consequently the linearization inherent to MSBDF2 or FEBDF2 leads to an inexact enforcement of the incompressibility condition. In other words, for the Liu model, the further $J$ deviates from 1 at a point, the worse the schemes presented in this paper (or any semi-implicit method that depends on similar linearization) are at enforcing the incompressibility condition within a neighborhood of that point. We therefore conclude that FEBDF2 is a robust and efficient integrator for incompressible nonlinear elastodynamics, particularly when paired with a quadratic formulation for the volumetric energy.

Using a nonlinear solver for the pressure equation would theoretically eliminate the problem of volume conservation by avoiding a linearization. A natural extension of the work in this paper is a comparison of the results obtained via FEBDF2 with an implicit solver that uses a Jacobian-Free Newton-Krylov (JFNK) scheme\cite{Chan1984}. JFNK schemes avoid much of the computational cost of typical nonlinear solvers by avoiding the explicit formation of a Jacobian matrix, and thus could potentially rival FEBDF2's efficiency when used to simulate large-scale dynamics in elastic solids while more accurately enforcing nonlinear incompressibility conditions.

\section{Acknowledgments}
We acknowledge the funding provided by the National Science Foundation via CSSI award NSF-OAC-1931516 and via the NSF-DMS-1929298 award to the Statistical and Applied Mathematical Sciences Institute. Numerical experiments were performed using computational facilities provided by the University of North Carolina at Chapel Hill through the Research Computing Division of UNC Information Technology Services.

\appendix
\section*{Appendices}
\addcontentsline{toc}{section}{Appendices}
\section{Truncation Error for MSBDF2}\label{appendix_A}
Because the pressure equation reduces to the enforcement of a constraint, it is not relevant in the following calculations. This would not be true if the incompressiblity constraint were enforced via a dynamic equation \cite{Rossi2016,Scovazzi2016}. For conciseness, functional dependence is only stated if it is relevant to the calculations. Upon substituting the true solution, $V^n = V(t^n)$, utilizing the truncated Taylor approximations
\begin{equation}
    \begin{aligned}
        {\textbf{V}}(t^n+\Delta t) \approx& {\textbf{V}}(t^n) + \Delta t \dot{\textbf{V}}(t^n) + \frac{\Delta t^2}{2}\ddot{\textbf{V}}(t^n) + \frac{\Delta t^3}{6}\dddot{\textbf{V}}(t^n)\\
        {\textbf{V}}(t^n-\Delta t) \approx& {\textbf{V}}(t^n) - \Delta t \dot{\textbf{V}}(t^n) + \frac{\Delta t^2}{2}\ddot{\textbf{V}}(t^n) - \frac{\Delta t^3}{6}\dddot{\textbf{V}}(t^n)\\
        \tilde{\mathbb{P}}^{\text{\text{dev}}}(t^n,t^n-\Delta t)\approx& \mathbb{P}^{\text{\text{dev}}}(t^n) + \Delta t \dot{\mathbb{P}}^{\text{\text{dev}}}(t^n) - \frac{\Delta t ^2}{2} \ddot{\mathbb{P}}^{\text{\text{dev}}}(t^n) \\
        p(t+\Delta t)\tilde{\mathbb{H}}(t^n,t^n-\Delta t)\approx& p(t^n)\mathbb{H}(t^n) 
        + \Delta t \frac{\partial }{\partial t}\left({p}(t^n)\mathbb{H}(t^n)\right)\\
        &+\Delta t^2\left[\dot{p}(t^n)\dot{\mathbb{H}}(t^n)+\frac{1}{2}(\ddot{p}(t^n)\mathbb{H}(t^n)-p(t^n)\ddot{\mathbb{H}}(t^n))\right]
    \end{aligned}
\end{equation}
in Equation (\ref{MSBDF2}a), and eliminating the time derivative of the conservation of linear momentum at $t=t^n$,
$$\rho_0 \ddot{\mathbf{V}}(t^n)=\nabla_0\cdot\left(\frac{\partial }{\partial t}{\mathbb{P}}(t^n)+\frac{\partial }{\partial t}(p(t^n)\mathbb{H}(t^n))\right) $$
the local truncation error at time $t^n$ can be expressed as
\begin{equation}
    \mathcal{L}^n = \Delta t^3\left[\frac{\rho_0}{12}\dddot{\textbf{V}}^n +\nabla_0\cdot\left(\ddot{\mathbb{P}}^n +\frac{1}{2}\ddot{p}^n\mathbb{H}^n-\dot{p}^n\dot{\mathbb{H}}^n-\frac{1}{2}p^n\ddot{\mathbb{H}}^n\right)\right] + O(\Delta t^4).
\end{equation}
Note that a lack of $\tilde{\mathbb{H}}(t^n,t^n-\Delta t)$ would eliminate cancellations in the error, reducing the scheme to first order.

\section{Truncation Error for FEBDF2}\label{appendix_B}
The work presented follows the same pattern established in Appendix \ref{appendix_A}. The Taylor approximations 
\begin{equation*}
    \begin{aligned}
        {\textbf{V}}(t^n+\Delta t) \approx& {\textbf{V}}(t^n) + \Delta t \dot{\textbf{V}}(t^n) + \frac{\Delta t^2}{2}\ddot{\textbf{V}}(t^n) + \frac{\Delta t^3}{6}\dddot{\textbf{V}}(t^n)\\
        {\textbf{V}}(t^n-\Delta t) \approx& {\textbf{V}}(t^n) - \Delta t \dot{\textbf{V}}(t^n) + \frac{\Delta t^2}{2}\ddot{\textbf{V}}(t^n) - \frac{\Delta t^3}{6}\dddot{\textbf{V}}(t^n)\\
        \tilde{\mathbb{P}}^{\text{\text{dev}}}(\mathbf{U}^*)\approx& \mathbb{P}^{\text{\text{dev}}}(t^n) + \Delta t \dot{\mathbb{P}}^{\text{\text{dev}}}(t^n) + \frac{\Delta t ^2}{2} \mathbf{V}(t^n)\cdot\partial_{\mathbf{UU}}{\mathbb{P}}^{\text{\text{dev}}}(t^n)\cdot\mathbf{V}(t^n)\\
        p(t+\Delta t)\tilde{\mathbb{H}}(\mathbf{U}^*)\approx& p(t^n)\mathbb{H}(t^n) 
        + \Delta t \frac{\partial }{\partial t}\left({p}(t^n)\mathbb{H}(t^n)\right)\\
        &+\Delta t^2\left[\dot{p}(t^n)\dot{\mathbb{H}}(t^n)+\frac{1}{2}(\ddot{p}(t^n)\mathbb{H}(t^n)+p(t^n)\mathbf{V}(t^n)\cdot\partial_{UU}{\mathbb{H}}(t^n)\cdot\mathbf{V}(t^n))\right]
    \end{aligned}
\end{equation*}
can be plugged into Equation (\ref{FEBDF2}). After eliminating Equation (\ref{CoLM}) and it's derivative in time, the local truncation error is found to be
\begin{equation}
    \mathcal{L}^n=\Delta t^3\left[\frac{\rho_0}{12}\dddot{\textbf{V}}^n -\nabla_0\cdot\left(\frac{1}{2} \mathbf{V}^n\cdot\partial_{\mathbf{UU}}\mathbb{P}^{\text{\text{dev}},n}\cdot\mathbf{V}^n +\frac{1}{2}\ddot{p}^n\mathbb{H}^n+\dot{p}^n\dot{\mathbb{H}}^n+\frac{1}{2}p^n\mathbf{V}^n\cdot\partial_{\mathbf{UU}}{\mathbb{H}}^{n}\cdot\mathbf{V}^n\right)\right] + O(\Delta t^4).
\end{equation}
This establishes FEBDF2 as a globally second order method.
 \bibliographystyle{ieeetr} 
 \bibliography{terrell_MSBDF2_FEBDF2}

\begin{thebibliography}{10}

\bibitem{Ita2022}
K.~Ita, M.~Silva, and R.~Bassey, ``Mechanical properties of the skin: {W}hat do we know?,'' {\em Current Cosmetic Science}, vol.~1, no.~1, pp.~70--76, 2022.

\bibitem{Mathur2001}
A.~Mathur, A.~Collinsworth, W.~Reichert, W.~Kraus, and G.~Truskey, ``Endothelial, cardiac muscle and skeletal muscle exhibit different viscous and elastic properties as determined by atomic force microscopy,'' {\em Journal of Biomechanics}, vol.~34, no.~12, pp.~1545--1553, 2001.

\bibitem{Schaefer2010}
R.~Shaefer, {\em Harris Shock and Vibration Handbook 6}, ch.~Mechanical Properties of Rubber.
\newblock McGraw-Hill, 2010.

\bibitem{Simo1984}
J.~Simo, R.~Taylor, and K.~Pister, ``Variational and projection methods for the volume constraint in finite deformation elasto-plasticity,'' {\em Computer Methods in Applied Mechanics and Engineering}, vol.~51, no.~1, pp.~177--208, 1984.

\bibitem{Coombs2010}
W.~Coombs and R.~Crouch, {\em 70-line 3{D} finite deformation elastoplastic finite-element code}, pp.~151--156.
\newblock Wiley, 2010.

\bibitem{Brink1996}
U.~Brink and E.~Stein, ``On some mixed finite element methods for incompressible and nearly incompressible finite elasticity,'' {\em Computational Mechanics}, vol.~19, no.~1, pp.~105--119, 1996.

\bibitem{Hughes2000}
T.~Hughes, {\em The Finite Element Method - {L}inear Static and Dynamic Finite Element Analysis}.
\newblock Dover, 2000.

\bibitem{Karabelas2020}
E.~Karabelas, G.~Haase, G.~Plank, and C.~Augustin, ``Versatile stabilized finite element formulations for nearly and fully incompressible solid mechanics,'' {\em Computational Mechanics}, vol.~65, no.~1, pp.~193--215, 2020.

\bibitem{Argyris1981}
J.~Argyris and S.~Symeonidis, ``Nonlinear finite element analysis of elastic systems under nonconservative loading-natural formulation. {P}art {I}. {Q}uasistatic problems,'' {\em Computer Methods in Applied Mechanics and Engineering}, vol.~26, no.~1, pp.~75--123, 1981.

\bibitem{Teran2005}
J.~Teran, E.~Sifakis, G.~Irving, and R.~Fedkiw, ``Robust quasistatic finite elements and flesh simulation,'' in {\em Proceedings of the 2005 ACM SIGGRAPH/Eurographics symposium on Computer animation}, 2005.

\bibitem{Quintal2011}
B.~Quintal, H.~Steeb, M.~Frehner, and S.~Schmalholz, ``Quasi-static finite element modeling of sieismic attenumation and dispersion due to wave-induced fluid flow in poroelastic media,'' {\em Journal of Geophysical Research: Solid Earth}, vol.~116, no.~B1, 2011.

\bibitem{Arriaga2007}
A.~Arriaga, J.~Lazkano, R.~Pagaldai, A.~Zaldua, R.~Hernandez, R.~Atxurra, and A.~Chrysostomou, ``Finite-element analysis of quasi-static characterization tests in thermoplastic materials: {E}xperimental and numerical analysis results correlation with {ANSYS},'' {\em Polymer Testing}, vol.~26, no.~3, pp.~284--305, 2007.

\bibitem{Leveque2007}
R.~LeVeque, {\em Finite Difference Methods for Ordinary and Partial Differential Equations: {S}teady-State and Time-Dependent Problems}.
\newblock SIAM, 2007.

\bibitem{Rossi2016}
S.~Rossi, N.~Abboud, and G.~Scovazzi, ``Implicit finite incompressible elastodynamics with linear finite elements: {A} stabilized method in rate form,'' {\em Computational Methods in Applied Mechanical Engineering}, vol.~311, pp.~208--249, 2016.

\bibitem{Kadapa2021}
C.~Kadapa, ``A novel semi-implicit scheme for elastodynamics and wave propagation in nearly and truly incompressible solids,'' {\em Acta Mechanica}, vol.~232, pp.~2135--2163, 2021.

\bibitem{Lahiri2005}
S.~Lahiri, J.~Bonet, J.~Peraire, and L.~Casals, ``A variationally consistent fractional time-step integration method for incompressible and nearly incompressible {L}agrangian dynamics,'' {\em International Journal for Numerical Methods in Engineering}, vol.~63, no.~10, pp.~1371--1395, 2005.

\bibitem{Verlet1967}
L.~Verlet, ``Computer "experiments" on classical fluids {I}. {T}hermodynamical properties of {L}ennard-{J}ones molecules,'' {\em Physical Review}, vol.~159, no.~1, pp.~98--103, 1967.

\bibitem{Gil2014}
A.~Gil, C.~Lee, J.~Bonet, and M.~Aguirre, ``A stabilized {P}etrov-{G}alerkin formulation for linear tetrahedral elements in compressible, nearly incompressible, and truly incompressible fast dynamics,'' {\em Computer Methods in Applied Mechanics and Engineering}, vol.~276, pp.~659--690, 2014.

\bibitem{Gottliev2003}
S.~Gottlieb and L.~Gottlieb, ``Strong stability preserving properties of {R}unge-{K}utta time discretization methods for linear constant coefficient operators,'' {\em Journal of Scientific Computing}, vol.~18, no.~1, pp.~83--109, 2003.

\bibitem{Scovazzi2016}
G.~Scovazzi, B.~Carnes, X.~Zeng, and S.~Rossi, ``A simple, stable, and accurate linear tetrahedral finite element for transient, nearly, and fully incompressible solid dynamics: {a} dynamic variational multiscale approach,'' {\em International Journal for Numerical Methods in Engineering}, vol.~106, no.~10, pp.~799--839, 2016.

\bibitem{Zeng2017}
X.~Zeng, G.~Scovazzi, N.~Abboud, C.~O., and S.~Rossi, ``A dynamic variational multiscale method for viscoelasticity using linear tetrahedral elements,'' {\em International Journal for Numerical Methods in Engineering}, vol.~112, no.~13, pp.~1951--2003, 2017.

\bibitem{Chung1994}
J.~Chung and J.~Lee, ``A new family of explicit time integration methods for linear and non-linear structural dynamics,'' {\em International Journal for Numerical Methods in Engineering}, vol.~37, no.~23, pp.~3961--3976, 1994.

\bibitem{Britz1997}
D.~Britz, ``Stability of the backward differentiation formula ({FIRM}) applied to electrochemical digital simulation,'' {\em Computers \& Chemistry}, vol.~21, no.~2, pp.~97--108, 1997.

\bibitem{Chow2021}
K.~Chow and S.~Ruuth, ``Linearly stabilized schemes for the time integration of stiff nonlinear {PDE}s,'' {\em Journal of Scientific Computing}, vol.~87, no.~3, p.~95, 2021.

\bibitem{Bashforth1883}
F.~Bashforth, {\em An Attempt to Test the Theories of Capillary Action by Comparing the Theoretical and Measured Forms of Drops of Fluid. With an Explanation of the Method of Integration Employed in Constructing the Tables Which Give the Theoretical Forms of Such Drops}.
\newblock Cambridge, 1883.

\bibitem{Keita2024}
S.~Keita, A.~Beljadid, and Y.~Bourgault, ``Implicit and semi-implicit second-order time stepping methods for the {R}ichards equation,'' {\em Advances in Water Resources}, vol.~148, p.~103841, 2021.

\bibitem{Ascher1995}
U.~Ascher, S.~Ruuth, and B.~Wetton, ``Implicit-explicit methods for time-dependent partial differential equations,'' {\em SIAM Journal on Numerical Analysis}, vol.~32, no.~3, pp.~797--823, 1995.

\bibitem{Varah1980}
J.~Varah, ``Stability restrictions on second order, three level finite difference schemes for parabolic equations,'' {\em SIAM Journal on Numerical Analysis}, vol.~17, no.~2, pp.~300--309, 1980.

\bibitem{Chorin1968}
A.~Chorin, ``Numerical solution of the {N}avier-{S}tokes equations,'' {\em Mathematics of Computation}, vol.~22, no.~104, pp.~745--762, 1968.

\bibitem{Tadmor2012}
E.~Tadmor, R.~Miller, and R.~Elliott, {\em Continuum Mechanics and Thermodynamics: {F}rom Fundamental Concepts to Governing Equations}.
\newblock Cambridge University Press, 2012.

\bibitem{Capecchi2007}
D.~Capecchi and G.~Ruta, ``{P}iola's contribution to continuum mechanics,'' {\em Archive for History of Exact Sciences}, vol.~61, no.~4, pp.~303--342, 2007.

\bibitem{Kim2011}
B.~Kim, S.~Lee, J.~Lee, S.~Cho, H.~Park, S.~Yeom, and S.~Park, ``A comparison among neo-{H}ookean model, {M}ooney-{R}ivlin model, and {O}gden model for chloroprene rubber,'' {\em International Journal of Precision Engineering and Manufacturing}, vol.~13, pp.~759--764, 2012.

\bibitem{Pence2015}
T.~Pence and K.~Gou, ``On compressible versions of the incompressible neo-{H}ookean material,'' {\em Mathematics and Mechanics of Solids}, vol.~20, no.~2, pp.~157--182, 2015.

\bibitem{Mooney}
M.~Mooney, ``A theory of large elastic deformation,'' {\em Journal of Applied Physics}, vol.~11, no.~9, pp.~582--592, 1940.

\bibitem{Ogden1972}
R.~Ogden, ``Large deformation isotropic elasticity -- on the correlation of theory and experiment for incompressible rubberlike solids,'' {\em Proceedings of the Royal Society of London. Series A, Mathematical and Physical Sciences}, vol.~326, p.~565–584, 1972.

\bibitem{Holzapfel2000}
G.~Holzapfel, T.~Gasser, and O.~R.W., ``A new constitutive framework for arterial wall mechancis and a comparative stud of material models,'' {\em Journal of Elasticity and the Physical Sciences of Solids}, vol.~61, no.~1, pp.~1--48, 2000.

\bibitem{Liu1992}
C.~Liu, G.~Hofstetter, and H.~Mang, ``3{D} finite element analysis of rubber-like materials at finite strains,'' {\em Engineering Computations}, vol.~11, no.~2, pp.~111--128, 1994.

\bibitem{Lvov2022}
V.~Lvov, F.~Senatov, A.~Veveris, V.~Skrybykina, and A.~Diaz~Lantada, ``Auxetic metamaterials for biomedical devices: {C}urrent situation, main challenges, and research trends,'' {\em Materials}, vol.~15, no.~4, p.~1439, 2022.

\bibitem{Courant1967}
R.~Courant, K.~Friedrichs, and H.~Lewy, ``On the partial difference equations of mathematical physics,'' {\em IBM Journal of Research and Development}, vol.~11, no.~2, pp.~215--234, 1967.

\bibitem{Xu2006}
C.~Xu and T.~Tang, ``Stability analysis of large time-stepping methods for epitaxial growth models,'' {\em SIAM Journal on Numerical Analysis}, vol.~44, no.~4, pp.~1759 -- 1779, 2006.

\bibitem{Brezzi1974}
F.~Brezzi, ``On the existence, uniqueness and approximation of saddle-point problems arising from {L}agrangian multipliers,'' {\em R.A.I.R.O. Analyse Numerique}, vol.~8, no.~R2, pp.~129--151, 1974.

\bibitem{Wang2020}
X.~Wang, Y.~Zou, and Q.~Zhai, ``An effective implementation for {S}tokes equation by the weak {G}alerkin finite element method,'' {\em Journal of Computational and Applied Mathematics}, vol.~370, 2020.

\bibitem{Xu2019}
J.~Xu, Y.~Li, S.~Wu, and A.~Bousquet, ``On the stability and accuracy of partially and fully implicit schemes for phase field modeling,'' {\em Computer Methods in Applied Mechanics and Engineering}, vol.~345, pp.~826--853, 2019.

\bibitem{Simo1992}
J.~C. Simo and N.~Tarnow, ``The discrete energy-momentum method. {C}onserving algorithms for nonlinear elastodynamics,'' {\em Zeitscrift fur Angewandte Mathematik und Physik}, vol.~43, no.~5, pp.~757 --792, 1992.

\bibitem{Kang2023}
Y.~Kang, ``An energy stable linear {BDF}2 scheme with variable time-steps for the molecular beam epitaxial model without slope selection,'' {\em Communications in Nonlinear Science and Numerical Simulation}, vol.~118, p.~107047, 2023.

\bibitem{Joly2007}
P.~Joly, {\em Numerical Methods for Elastic Wave Propagation}, ch.~6.
\newblock Springer, 2007.

\bibitem{Voet2023}
Y.~Voet, E.~Sande, and A.~Buffa, ``A mathematical theory for mass lumping and its generalization with applications to isogeometric analysis,'' {\em Computer Methods in Applied Mechanics and Engineering}, vol.~410, p.~116033, 2023.

\bibitem{Taylor1973}
C.~Taylor and P.~Hood, ``A numerical solution of the {N}avier-{S}tokes equations using the finite element method,'' {\em Computers \& Fluids}, vol.~1, no.~1, pp.~73--100, 1973.

\bibitem{Brezzi1992}
F.~Brezzi, M.~Bristau, L.~Franca, M.~Mallet, and G.~Rog\'e, ``A relationship between stabilized finite element methods and the {G}alerkin method with bubble functions,'' {\em Computer Methods in Applied Mechanics and Engineering}, vol.~98, no.~1, pp.~117--129, 1992.

\bibitem{Cohen2001}
G.~Cohen, P.~Joly, J.~Roberts, and N.~Tordjman, ``Higher order triangular finite elements with mass lumping for the wave equation,'' {\em SIAM Journal on Numerical Analysis}, vol.~28, no.~6, pp.~2047--2078, 2001.

\bibitem{Hestenes1952}
M.~Hestenes and E.~Stiefel, {\em Methods of conjugate gradients for solving linear systems}.
\newblock National Bureau of Standards, 1952.

\bibitem{Falgout2006}
R.~Falgout, ``An introduction to algebraic multigrid,'' {\em Computing in Science \& Engineering}, vol.~8, no.~6, pp.~24--33, 2006.

\bibitem{dealii}
P.~C. Africa, D.~Arndt, W.~Bangerth, B.~Blais, M.~Fehling, R.~Gassmöller, T.~Heister, L.~Heltai, S.~Kinnewig, M.~Kronbichler, M.~Maier, P.~Munch, M.~Schreter-Fleischhacker, J.~P. Thiele, B.~Turcksin, D.~Wells, and V.~Yushutin, ``The deal.{II} library, version 9.6,'' {\em Journal of Numerical Mathematics}, vol.~32, no.~4, pp.~369--380, 2024.

\bibitem{paraview}
\url{https://www.paraview.org}.

\bibitem{matlab}
\url{https://www.mathworks.com/products/matlab.html}.

\bibitem{Chan1984}
T.~Chan and K.~Jackson, ``Nonlinearly preconditioned {K}rylov subspace methods for discrete {N}ewton algorithms,'' {\em SIAM Journal on Scientific and Statistical Computing}, vol.~5, no.~3, pp.~533--542, 1984.

\end{thebibliography}

\end{document}